\documentclass[12pt,a4paper]{article}
\usepackage{times,amsmath,amsthm,amssymb,amsfonts,tabularx,hhline,rotating,enumerate,comment}
\usepackage{graphicx}
\usepackage{hyperref}
\usepackage[usenames]{color}
\newcommand{\SIG}{\Sigma}
\newcommand{\GG}{\Gamma}
\newcommand{\gp}{geodesically perfect}

\newcommand{\IFF}{if and only if } 
\newcommand{\set}[2]{\left\{\, \mathinner{#1}\vphantom{#2}\; \left|\; \vphantom{#1}\mathinner{#2} \right.\,\right\}}
\newcommand{\oneset}[1]{\left\{\, \mathinner{#1} \,\right\}}
\newcommand{\smallset}[1]{\left\{\mathinner{#1}\right\}}
\newcommand{\abs}[1]{\left|\mathinner{#1}\right|}

\newcommand{\gen}[1]{\left< \mathinner{#1} \right>}

\newcommand\IRR{\mathop\mathrm{IRR}}
\newcommand\eps{\varepsilon}

\newcommand\ov[1]{\overline{#1}}

\newcommand{\Oh}{\mathcal{O}}

\newcommand\ZZ{\mathbb{Z}}
\newcommand\NN{\mathbb{N}}

\newcommand\HNN{\mathrm{HNN}}

\newcommand\RAS[2]{\overset{#1}{\underset{#2}{\Longrightarrow}}}
\newcommand\ras[2]{\overset{#1}{\underset{#2}{\longrightarrow}}}
\newcommand\LAS[2]{\overset{#1}{\underset{#2}{\Longleftarrow}}}
\newcommand\las[2]{\overset{#1}{\underset{#2}{\longleftarrow}}}
\newcommand\DAS[2]{\overset{#1}{\underset{#2}{\Longleftrightarrow}}}
\newcommand\das[2]{\overset{#1}{\underset{#2}{\longleftrightarrow}}}
\newcommand\OUTS[5]{#1
\overset{#2}{\underset{#3}{\Longleftarrow}} #4
\overset{#2}{\underset{#3}{\Longrightarrow}} #5}
\newcommand\INS[5]{#1
\overset{#2}{\underset{#3}{\Longrightarrow}} #4
\overset{#2}{\underset{#3}{\Longleftarrow}} #5}
\newcommand\RA[1]{{\underset{#1}{\Longrightarrow}}}
\newcommand\LA[1]{{\underset{#1}{\Longleftarrow}}}
\newcommand\DA[1]{{\underset{#1}{\Longleftrightarrow}}}
\newcommand\OUT[4]{#1
{\underset{#2}{\Longleftarrow}} #3
{\underset{#2}{\Longrightarrow}} #4}
\newcommand\IN[4]{#1
{\underset{#2}{\Longrightarrow}} #3
{\underset{#2}{\Longleftarrow}} #4}

\newtheorem{theorem}{{\bf Theorem}}[section]
\newtheorem{corollary}[theorem]{{\bf Corollary}}
\newtheorem{definition}[theorem]{{\bf Definition}}
\newtheorem{example}[theorem]{{\bf Example}}
\newtheorem{lemma}[theorem]{{\bf Lemma}}
\newtheorem{proposition}[theorem]{{\bf Proposition}}
\newtheorem{remark}[theorem]{{\bf Remark}}
\newtheorem{problem}[theorem]{{\bf Problem}}
\newtheorem{conjecture}[theorem]{{\bf Conjecture}}

\newenvironment{am}{\noindent\color{blue} AM: }{}
\newenvironment{ajd}{\noindent\color{red} AJD }{}
\newenvironment{vd}{\noindent\color{green} VD }{}

\newcommand{\ad}[1]{ \begin{ajd} #1 \end{ajd}}

\newcommand{\be}{\begin{enumerate}}
\newcommand{\ee}{\end{enumerate}}
\begin{document}

\title{Geodesic rewriting systems and pregroups\thanks{
Part of this work was begun in 2007 when the first and third author where
at the CRM (Centro Recherche Matem{\`a}tica, Barcelona) at the  invitation of
 Enric Ventura.}}
\author{Volker Diekert${}^1$ \and Andrew J.~Duncan${}^2$ \and Alexei Miasnikov${}^3$\\
${}^1$Universit\"at Stuttgart,
Universit\"atsstr. 38\\ D-70569 Stuttgart, Germany\\
${}^2$Newcastle University, Newcastle upon Tyne\\ NE1 7RU, United Kingdom\\
${}^3$McGill University,
Montreal,
Canada,
H3A 2K6}

\date{\today}

\maketitle

\begin{abstract}
In this paper we study rewriting systems for groups and monoids,
focusing on situations where finite convergent systems may be
difficult to find or do not exist. We consider systems which
have no length increasing rules and are confluent and then systems
in which the length reducing rules lead to geodesics. Combining
these properties we arrive at our main object of study which
we call geodesically perfect rewriting systems. We show that
these are well-behaved and convenient to use, and give several
examples of classes of groups for which they can be constructed
from natural presentations. We describe a Knuth-Bendix completion
process to construct such systems, show how they may
be found with the help of Stallings' pregroups and conversely may
be used to construct such pregroups.
\end{abstract}

\tableofcontents

\section{Introduction}\label{sec:intro}
A presentation of a group or monoid may be thought of as a rewriting
system which, in certain cases may give rise to algorithms for solving 
classical algorithmic problems.
For example if the rewriting system is finite and convergent (that
is confluent and terminating)
then it can be used to solve the word problem and
to find normal forms for elements of the group.
This is one reason for the importance of convergent rewriting systems
in group theory. However there are many groups for which the natural
presentations do not give rise to convergent rewriting systems, but which
are none the less well behaved, algorithmically tractable groups. In this
paper we investigate properties of rewriting systems, which are not
in general finite or terminating, but which all the same give algorithms
for such tasks as solving the 
word problem, computation of normal forms or computation of geodesic
representatives of group elements. We contend  that the resulting
algorithms are often more convenient and practical than those arising
from more conventional finite convergent systems.

Rewriting methods in algebra have a very long and rich history.
In groups and semigroups they are usually related to the word
problem and take their roots in the ground breaking works of
Dehn and Thue (not to mention the classical Euclidean and
Gaussian elimination algorithms!). Several famous algorithms in group
theory are in fact  particular types of string rewriting processes:
the Nielsen method  in  free groups, Hall collection in nilpotent
and polycyclic groups, the Dehn algorithm in small cancellation and
hyperbolic groups, Tits rewriting in Coxeter groups, convergent
rewriting systems for finite groups, and so on.  In rings and algebras
rewriting methods appear as a major tool in computing normal forms
of elements \cite{Newman42,Shirshov62,Buchberger65}, in solving
the word and ideal-membership problems. These techniques emerged
up independently in various branches of algebra at different times
and  under different names  (the diamond lemma, Gr\"obner or Shirshov
bases, Buchberger's algorithm and $S$-polynomials, for instance).  They
have gained
prominence with the progress of practical computing,  as
real applications have become available.  Notably,  crucial  developments
in methods  of computational algebra  originated in commutative
algebra and algebraic geometry, with  Buchberger's  celebrated
algorithm
and related computational techniques, which revolutionised the whole area of
applications. We refer to \cite{CoxLittleOshea07}, and the
references therein, for more details.

From the theoretical view point the main shift in the
paradigm came with the seminal paper of  Knuth and
Bendix \cite{KnuthBendix70}.
In this paper they
introduced a process, now known as the  Knuth-Bendix (KB) procedure,
which
unified the field of rewriting techniques in 
(universal) algebra.  This KB procedure gives a solid theoretical basis for practical implementations, even the procedure itself may lead to
non-optimal algorithms for solving word problems.

Roughly speaking
a KB procedure takes as input a finite  system of identities 
(between terms) and a computable (term) ordering such that the identities can 
be read as a finite set of directed rewrite rules. 
Using the crucial concept of \emph{critical pairs} the 
procedure adds in each round more and more rules, and it stops only  if the 
system is \emph{completed}. Thus the KB procedure 
attempts to construct an equivalent convergent (term) rewriting system: which
in particular allows unique normal forms to be found by a simple strategy. In case of termination  we obtain a solvable word problem. 

In the case of commutative algebra this concept 
can be viewed as  Buchberger's algorithm and termination is guaranteed.
In case  of algebraic structures like
 groups or monoids we have a special case of term rewriting systems
 since the rewriting process is based on strings. 
 (Formally, monoid generators are read as unary function symbols,
 and the neutral element is read a constant.)

As has been mentioned above, the history of rewriting systems in monoids
and groups is about one hundred years old, with the main focus on
convergent rewriting systems and algorithms for computing normal forms.
Any presentation $M = \langle \Gamma \mid \ell_i = r_i (i \in I)\rangle$
of a monoid $M$ gives a rewriting system $S = \{ \ell_i \to r_i, i \in I\}$
which defines
$M$ via the congruence relation it generates on the free monoid $\Gamma^\ast$.
Every rule $\ell \to r \in S$ allows one to rewrite a word $u\ell v$ into the word $urv$ 
and this gives a (non-deterministic) word rewriting procedure
associated with $S$.
If the system $S$ is convergent (see Section \ref{sec:rewrite}) then
this rewriting system describes
a  deterministic algorithm which computes the normal forms of elements,
thus solving the word problem in the monoid $M$.
This yields the major interest in finite convergent  systems.  Many groups
are known to allow  finite convergent systems (for example Coxeter groups,
polycyclic groups, some small cancellation groups: see books
\cite{HoltEickObrien05,Sims94,Chenadec86} for more examples and details).
The primary task here is to find a finite convergent
system for a given finitely presented monoid, assuming that such a
system exists. In principle, the KB procedure performs this task. However,
several obstacles may present themselves. By design, to start the KB procedure
one has to fix in advance an ordering on $\Gamma^\ast$, with particular
properties, as described in Section \ref{sec:convergent}. This
may seem like a minor hurdle, but the difficulty is that, even for well understood
groups, with two  orderings which look very much alike, it may happen
that using the first the
KB process halts and outputs a convergent system while with respect to the
the second there exists no finite convergent system:
see Example \ref{ex:badorders} below.
Furthermore, the existence of a finite convergent system also depends
on a choice of the set of generators of the group.
This means that for KB to succeed one has to make a the right
choice of a set of generators $\Gamma$ and of an ordering on
$\Gamma^\ast$. 
In fact \cite{odunlaing83}
in general the
problem of whether or not a given finitely presented group can
be defined by a finite convergent rewriting system  
is undecidable. In addition,
even when restricted to instances where 
the generators and the order have been chosen so that
the KB process will halt giving  a finite convergent
rewriting system, there may be no be effectively computable upper bound on the
running time of the  KB procedure. 
To make things even more interesting, having a finite convergent
rewriting system $S$ does not guarantee  a fast solution of the word
problem in the monoid $M$ (see Section \ref{sec:convergent}).
All these results show that the KB  process for finite convergent systems,
while being an important theoretical tool, is not a
panacea for problems in computational algebra.

As a first step towards resolving some of these difficulties
we consider, in Section
\ref{sec:preperfect},
the class of preperfect
rewriting systems: that is those which are confluent and
have no length increasing
rules.  These restrictions are enough to allow solution
of the word problem and to find geodesic representatives and, as examples
show, such systems are common in geometric group theory.
In fact in Section \ref{sec:examples} we describe preperfect rewriting systems for Coxeter groups, graph groups,
HNN-extensions and free products with amalgamation.
One disadvantage of these systems is that it is undecidable whether are a
finite rewriting system is preperfect or not \cite{NarendranMcNaughton84}
(see Theorem \ref{thm:prepundec}).

Another desirable property of rewriting systems is that they
should be geodesic; meaning that shortest representatives of
elements can be found by applying
only the length reducing rules of the system. A group defined by  a
finite geodesic rewriting system has solvable word problem and
in  \cite{GilHHR07} these groups are characterised as
the finitely generated virtually free groups.
However, as we show in Section \ref{sec:geodesic}, the question
of a whether or not a finite rewriting system is geodesic is undecidable.

Combining properties of preperfect and geodesic rewriting systems
we arrive at geodesically perfect rewriting systems
(defined in Section \ref{sec:geodesicperf}). These were first investigated
by Nivat and Benois \cite{nibe71} where they were  called
\emph{quasi-parfaites}.
Elsewhere these rewriting systems are also known as
\emph{almost confluent}, see e.g. \cite{bo93springer} but
here we prefer the notation  geodesically perfect since
these systems are designed to deal with geodesics in groups and monoids.
In  \cite{nibe71},
it was shown that
the property of being
geodesically perfect
is decidable for finite
systems. This leads to a new Knuth-Bendix completion procedure
for constructing
geodesically perfect systems as we  explain below. One advantage
of this KB process is that it requires no choice of ordering, using
only the partial order given by word length in $\Gamma^\ast$.

Among the examples of Section \ref{sec:examples} are rewriting  systems for
amalgamated products and HNN-extensions. As several
several important  frameworks have been developed to unify the studies of
such groups (Bass-Serre Theory, pregroups
and relatively hyperbolic groups,for example)
it is natural to look for a unified theory of
rewriting systems covering HNN-extensions and amalgamated products.
In this paper, following Stallings \cite{Stallings71,Stallings87}, we approach this unification question
from a combinatorial view-point via pregroups and their universal groups:
which seem to lend themselves naturally to algorithmic and model theoretic
problems.
 Intuitively, a pregroup can be viewed as a ``partial group'', that is,
a set $P$ with a partial (not everywhere defined) multiplication
$m:P \times P \to P$, or a piece of the multiplication table of some group,
that satisfies some particular axioms. In this case the universal
group $U(P)$ can be described as the group defined by the presentation
with a generating set $P$ and a set of relations $m(x,y) = z$ for all
$x,y \in P$ such that $m(x,y)$ is defined and equal to $z$.
On the other hand, Stallings proved that  $U(P)$ can be realized
constructively as the set of all $P$-reduced forms
(reduced sequences of elements of $P$)  modulo a suitable
equivalence relation and a naturally defined multiplication.
We discuss these definitions in detail in Section \ref{sec:pregroup}.

In Section \ref{sec:rewrite-U(P)}  we show how the existence of
a pregroup allows us to construct a preperfect rewriting system
for the universal group.  Moreover, we show in Theorem
\ref{th:preperfect-pregroups} that this system is geodesically perfect.
In this way   pregroups may play a
role in clarifying   completion procedures of KB type. In particular,
completing a given presentation  (in terms of generators and relators)
of a group $G$ to a larger   presentation,  which is  a pregroup,
amounts to a construction of a  geodesically
perfect rewriting system for $G$.

As an application of these results we obtain a slight strengthening
of the result of \cite{GilHHR07}.
It is known that a group $G$ is virtually free if and only if   $G = U(P)$ for a finite pregroup $P$ \cite{Rimlinger87b} and combining this result with
Theorem \ref{th:preperfect-pregroups} we see that a group is
finitely generated, virtually free if and only if it is defined
by a geodesically perfect rewriting system (Corollary \ref{virtfreegp}).

\section{Rewriting techniques}\label{sec:rewrite}
\subsection{Basics}

In this section we recall the basic concepts from string rewriting.
We use rewriting techniques as a tool to prove that certain constructions
have the expected properties.

A \emph{rewriting relation} over a set $X$ is a binary relation
$\RA{}\subseteq X \times X$. We denote by $\RAS *{}$ the
reflexive and transitive closure of $\RA{}$, by
$\DA{}$ its symmetric closure and by  $\DAS*{}$ its
symmetric,
reflexive, and transitive closure. We also write
$y \LA{}x$ whenever $x \RA{} y$, and we write $x \RAS{\leq k} {} y$
whenever we can reach $y$ in at most $k$ steps from $x$.

\begin{definition}
The  relation
$\RA{}\subseteq X \times X$ is called:
\begin{enumerate}[i)]
\item \emph{strongly confluent}, if $\OUT y{}xz$ implies
$\INS y{\leq 1}{}wz$ for some $w$,
\item \emph{confluent}, if $\OUTS y*{}xz$ implies
$\INS y{*}{}wz$ for some $w$,
\item \emph{Church-Rosser}, if $y\DAS * {}z$ implies
$\INS y{*}{}wz$ for some $w$,
\item \emph{locally confluent}, if $\OUT y{}xz$ implies
$\INS y{*}{}wz$ for some $w$,

\end{enumerate}
\end{definition}

The following facts are well-known  and
can be found in  several text
books (see for example,  \cite{bo93springer,jan88eatcs}).

\begin{enumerate}[1)]
\item Strong confluence  implies confluence.
\item Confluence  is equivalent to Church-Rosser.
\item Confluence implies local confluence, but the converse is false,
in general.
\end{enumerate}

\subsection{Rewriting in monoids}

Rewriting systems over monoids (and in particular over groups)
play an important  part  in algebra.
Let $M$  be a  monoid.  A \emph{rewriting system
} over $M$ is a binary relation
$S \subseteq M \times M$. It defines the rewriting relation
$\RA S \subseteq M \times M$ such that
\begin{align*}
x   \RA S y  \;\mbox{if and only if } \; x=p\ell q, \; y= prq \; \mbox{ for some } \; (\ell,r)
\in S.
\end{align*}
The relation  $\DAS * S \subseteq M \times M$ is a congruence on $M$, hence
the quotient set $M/\DAS * S $ forms a monoid with respect to the multiplication induced from $M$. We denote it by $M/\set{\ell=r}{(\ell,r)\in S} $ or, simply by  $M/S $. Two rewriting systems $S$ and $T$ over a monoid $M$ are termed  {\em equivalent} if $\DAS * S  = \DAS * T$, i.e., $M_S = M_T$.

We say that a rewriting system $S$ is strongly confluent (or confluent, etc) if
the relation $\RA S$ has the corresponding property. Instead of $(\ell,r) \in S$ we
also write $\ell {\underset{}{\longrightarrow}} r \in S$ and $\ell
{\underset{}{\longleftrightarrow}} r \in S$ in order to indicate that
both $(\ell,r)$ and $(r,\ell)$ are in $S$.

We say that a word $w$ is $S$-\emph{irreducible} (sometimes we omit $S$  here), if no left-hand side
$\ell$
of $S$ occurs in $w$ as a factor. Thus, if $w$ is irreducible, then
$w \RAS*S w'$
implies $w = w'$. The set of all irreducible words is denoted
by $\IRR(S)$.

In order to compute with monoids (in particular, groups) we usually specify
a choice of  monoid generators $\GG$, sometimes called an \emph{alphabet}.
For groups we often assume that $\GG$ is closed under inversion, so  $\GG = \Sigma \cup \Sigma^{-1}$ where
$\Sigma$ is a set of group generators. For an alphabet   $\GG$
we denote by $\Gamma^*$ the free monoid with basis $\Gamma$.
Throughout, $1$ denotes the neutral element in monoids or groups.
In particular, $1$ is also used to denote the  \emph{empty word}
in a free monoid $\Gamma^*$. If we can write 
$w= xuy$, then we say that $u$ is a factor of $w$. For free monoids a
factor is sometimes also called a \emph{subword}, 
but this might lead to confusion because other authors understand
b a subword simply a subsequence or \emph{scattered subword}. 

Rewriting systems $S$ over a free monoid $\Gamma^*$  are sometimes called \emph{string rewriting systems} or \emph{semi-Thue systems}.
In this  case the quotient $\Gamma^*/S$ has the standard monoid presentation $\langle \Gamma \mid \set{\ell=r}{ (\ell,r)\in S} \rangle$.
We say that  a string rewriting system $S$ defines a monoid $M$ if $\Gamma^*/S$  is isomorphic to $M$. In addition, if  $P$ is a property of  rewriting systems (Church-Rosser, strongly confluent, confluent,  etc.) we say that a monoid $M$ has a $P$-presentation if it can be defined by a system with property $P$.

For groups two types of presentations via generators and relators arise:
monoid presentations, described above, and group presentations,
typical in combinatorial group theory and topology.  More precisely, we say that $G = \Gamma^*/S$ is  a monoid presentation of a group $G$ if the alphabet $\Gamma$ is of the form $\GG = \Sigma \cup \Sigma^{-1}$, where
$\Sigma$ is a set of group generators, and $\Sigma^{-1}  = \{ \sigma^{-1} \mid \sigma \in \Sigma\}$ is the set of formal inverses of $\Sigma$ (in which case $\Gamma^\ast$ is a the free monoid  with an involution $\sigma \to \sigma^{-1}$).
 Given a group  presentation $\langle X \mid R \rangle$ of a group $G$  one can easily obtain a monoid presentation of $G$ by adding the formal inverses $X^{-1}$ to the set of generators $X$ of $G$ and the ``trivial'' relations $xx^{-1} = 1, x^{-1}x = 1, x \in X\}$ to the relators of $G$.
We consider here monoid presentations of groups, except where explicitly
indicated otherwise.

\subsection{Convergent rewriting systems}\label{sec:convergent}

In this section we briefly discuss {\em convergent} (or {\em complete}) rewriting systems, which play an important role in algebra due to their relation to normal forms.

A relation
$\RA{}\subseteq X \times X$ is called  \emph{terminating} (or \emph{Noetherian}), if every infinite chain
\begin{align*}
  x_0 \RAS *{} x_1 \RAS *{} \cdots x_{i-1} \RAS *{} x_i  \RAS *{} \cdots
\end{align*}
becomes stationary.

There are two typical sources of terminating string rewriting systems $S  \subseteq \Gamma^* \times \Gamma^*$. Systems of the first type are {\em length-reducing}, i.e., for any rule $\ell \to r \in S$ one has $|\ell| >|r|$, where $|x|$ is the length of a word $x \in \Gamma^*$. Systems of the second type are {\em compatible} with  a given {\em reduction}  ordering $\succ$ on $\Gamma^*$, which means that if $\ell \to r \in S$ then $\ell \succ r$. Recall that a reduction ordering on $\Gamma^*$  is a well-ordering preserving left and right multiplication (i.e. if $u \succ v$ then $aub \succ avb$ for any $a,b \in \Gamma^*$). Clearly, such systems are terminating.  In fact, the condition that $S$ is compatible with  some partial order, $\succ$, 
preserving left and right multiplication
is just a reformulation of the terminating property. Indeed, if $S$ is terminating then there is a binary relation $\succ_S$ on $\Gamma^*$ defined by  $u \succ_S v$ if and only if $u  \RAS *{S}  v$. In this case  $\succ_S$ is a partial well-founded ordering (no infinite descending chains),  such that  $\ell \succ_S r$ for any rule $\ell \to r \in S$. Moreover, the converse is also true. (The condition that $\succ$ is total is not
needed here but is required in running  the Knuth-Bendix completion procedure, see below).

 A relation $\RA{}\subseteq X \times X$ is called \emph{convergent} (or \emph{complete})  if it is locally  confluent and terminating.
The following properties are crucial. Let $S$
be a convergent rewriting system.
\begin{enumerate} [1)]
 \item $S$ is confluent  (see for example,  \cite{bo93springer,jan88eatcs}).

 \item Every $\DAS*{S}$ equivalence class in $\Gamma^*$ contains a unique $S$-reduced word (a word to which no rule from $S$ is applicable).

 \item If $S$ is finite then for a given word $w \in \Gamma^*$ one can effectively find its unique $S$-reduced form (just by subsequently rewriting the word $w$ until the result is $S$-reduced).
\end{enumerate}

 The results above show that if a monoid $M$ has a finite convergent presentation then the word problem in $M$ as well as the problem of finding the normal forms, is decidable. This explains popularity of convergent systems in algebra. There are many examples of groups that have finite convergent presentations:
finite groups, polycyclic group,  free groups, some geometric groups  (see \cite{Sims94, ech92, Chenadec86} for details)

One of the major results on convergent systems  concerns the Knuth-Bendix procedure (KB) (see \cite{bo93springer} for general rewriting systems and \cite{Sims94,ech92} for groups), which can be stated as follows.  Let $\succ$ be a
reduction  well-ordering on $\Gamma^*$  and
$S \subseteq \Gamma^* \times \Gamma^*$ a finite rewriting system
compatible with  $\succ$. If there exists a finite convergent rewriting system $T \subseteq \Gamma^* \times \Gamma^*$ compatible with  $\succ$ which is equivalent to $S$, then, in finitely many steps, the Knuth-Bendix procedure KB finds a finite convergent rewriting system $S^\prime \subseteq \Gamma^* \times \Gamma^*$ compatible with  $\succ$ which is also equivalent to $S$.

There are three  principle remarks due  here.

\begin{remark}
The time complexity of the word problem in a monoid $M_S$ defined by a finite convergent system $S$ may be of an arbitrarily high
complexity
\cite{OttoK97}.

\end{remark}

\begin{remark}
 	It may happen that the word problem in a monoid $M_S$ defined by a finite convergent system $S$ is decidable in polynomial time, whereas
the complexity of the standard rewriting algorithm that finds  the $S$-reduced forms of words can be of an arbitrarily high  complexity \cite{OttoK97}.
\end{remark}

These  remarks show that convergent rewriting systems may not be the best tool to deal with complexity issues related to the word problems and normal forms in monoids.

\begin{remark}
 The Knuth-Bendix procedure really depends on the chosen ordering $\succ$. The following example shows that in a free Abelian group of rank two the KB procedure relative to one length-lexicographic ordering results in a finite convergent presentation, while another length-lexicographic ordering does not allow any finite convergent presentations for the same group.
 \end{remark}

\begin{example}[\cite{Epstein}, page 127]\label{ex:badorders} Let $G$ be the free Abelian group  given by the following monoid presentation.
$$\langle x,y,x^{-1}, y^{-1} \mid xy = yx , xx^{-1} = x^{-1}x = yy^{-1} = y^{-1}y = 1\rangle.$$
Then the KB procedure with respect to the length-lexicographic ordering induced by the ordering $x < x^{-1} < y < y^{-1}$ of the generators outputs a finite convergent system defining $G$:
$$xx^{-1} \RA{} 1, x^{-1}x \RA{} 1,  yy^{-1} \RA{} 1,  y^{-1}y \RA{} 1, $$
$$yx \RA{} xy, y^{-1}x \RA{} xy^{-1}, yx^{-1} \RA{} x^{-1}y, y^{-1}x^{-1} \RA{} x^{-1}y^{-1}.$$
However,  there are no finite convergent systems defining $G$ and compatible with the  length-lexicographic ordering $x <y < x^{-1} <y^{-1}$.
 \end{example}

 Therefore, even if a finite convergent presentation for a monoid $M$ exists it might be hard to find it using the Knuth-Bendix procedure.
In addition {\'O}'D{\'u}n\-laing \cite{odunlaing83} has shown that
the
problem of whether or not a given finitely presented group can
be defined by a finite convergent rewriting system  
is undecidable.

It is not hard to see that all finitely generated commutative monoids
have a finite convergent presentation, \cite{die86tcsb}. However,
this is demands enough generators, in general. For example,
a free Abelian groups of rank $k$ can be generated as a monoid by an 
alphabet of size $k+1$, but in  order to find a 
finite convergent system for it we need at least $2k$ generators, 
see \cite{die86tcs}. Another nice example of this kind 
is the non-commutative semi-direct product of $\ZZ$ by  $\ZZ$. 
Even as  a monoid we need just two generators $a$ and $b$ and 
one relation $abba = 1$. There is no finite convergent system 
$S \subseteq \oneset{a,b}^* \times  \oneset{a,b}^*$ such that
$\oneset{a,b}^*/\oneset{abba= 1} =  \oneset{a,b}^*/S$, but clearly such 
systems exists if we spend more generators. See \cite{jan88eatcs} for more details about this example. 

We finish the section with a few open problems.

\begin{problem}
 Is it true that every hyperbolic group has a finite convergent presentation?
\end{problem}
 It is known that some hyperbolic groups have finite convergent presentations, for example, surface groups \cite{Chenadec86}.

\begin{problem}
 Is it true that every finitely generated fully residually free  group has a finite convergent presentation?
\end{problem}

The next two problems are from \cite{OttoK97}.

\begin{problem}
Do all automatic groups have finite convergent presentations?
\end{problem}

\begin{problem}
Do all one-relator groups have finite convergent presentations?
\end{problem}

Notice that all the groups above satisfy the homological condition
$FP_\infty$; which is the main   known  condition necessary
for a group  to have a finite convergent presentation, see \cite{SquOttKob94, Squier1987a}.

\subsection{Computing with infinite systems}
\label{subsec:computing}

In this section we discuss computing with  infinite systems. 
 An infinite string rewriting system  $S \subseteq \Gamma^* \times \Gamma^*$ can be used in computation if it satisfies some  natural conditions.
Firstly, one has to be able to recognise if a given pair $(u,v) \in \Gamma^* \times \Gamma^*$ gives a rule $u \to v \in S$ or not, i.e., the system $S$ must be  a recursive subset of $\Gamma^* \times \Gamma^*$. We call such  systems  {\em recursive}. Secondly, to rewrite with $S$ one has to be able to check if for a given $u \in \Gamma^*$ there is a rule $\ell \to r \in S$ with $\ell = u$, so we assume that the set $L(S)$ of the left-hand sides of the rules in $S$ is a recursive subset of $\Gamma^*$. Systems satisfying these two conditions are termed {\em effective} rewriting systems. Clearly, every finite system is effective. Notice also, that every recursive {\em non-length-increasing} system $S$ (i.e., $|\ell| \geq |r|$ for every rule $\ell \to r \in S$) is effective. Indeed, given $u \in \Gamma^*$ one can check if a rule $u \to v$ is in $S$ or not for all words $v$ with $|v| \leq |u|$, thus effectively verifying whether $u \in L(S)$ or not. 

The argument above shows that for a recursive non-length-increasing system $S$ one can effectively enumerate all the rules in $S$ in such a way
 \begin{equation}
 \label{eq:enumeration}
  \ell_0 \to r_0, \ell_1 \to r_1, \ldots ,\ell_i \to r_i, \ldots
  \end{equation}
   that if $i < j$ then $\ell_i \preceq \ell_j$ in the
length-lexicographical
ordering $\preceq$ and also if $\ell_i = \ell_j$ then $r_i \preceq r_j$. We call this enumeration of $S$ {\em standard}. 

\begin{proposition}
\label{pr:eff}
Let $S$ be an infinite  effective convergent system. Then the word problem in the monoid $M_S$ defined by $S$ is decidable.
 \end{proposition}
\begin{proof}
Given a word $u \in \Gamma^*$ one can start the rewriting process applying rules from $S$. Indeed, for a given  factor $w$ of $u$ one can check if $w \in L(S)$ or not, thus enumerating all factors of $u$ one can either find a factor $w$ of $u$  with $w \in L(S)$ or prove that $u$ is $S$-irreducible. If such $w$ exists one can enumerate  all pairs $(w,v)$ with $v \in \Gamma^*$  and check one by one if $(w,v) \in S$ or not. This procedure eventually terminates with a rule   $w \to v  \in S$. Now one can apply this rule to $u$ and rewrite $u$ into $u_1$.    Applying again this procedure to $u_1$ one eventually arrives at a 
unique $S$-irreducible word $\widehat{u}$.  To check if two words are equal in the monoid $M_S$ one can find their $S$-irreducibles and check whether they are equal 
or not.

\end{proof}

There are various modifications of the  algorithm described above, that work for other types of, not necessarily convergent, 
infinite systems. We consider some of these below.

\section{Length-reducing and Dehn systems}
\subsection{Finite length-reducing systems}

In this section we study a very particular type of rewriting system,
called {\em length-reducing} systems, where, for every rule $\ell \to r$
one has $|\ell| > |r|$.  The main interest in length-reducing systems  comes from the fact that, contrary to the case of finite convergent systems, the algorithm for computing the reduced forms is fast. 

\begin{lemma} \cite{boo82} If $S$ is a finite length-reducing string rewriting system, then
irreducible descendants of a given word can be computed in linear time
(in the length of the word).
\end{lemma}
This result is well-known,  we use it in many parts of the paper,  and it can be seen easily as follows.
\begin{proof} First, we choose some
$\eps > 0$ such that $(1-\eps) |\ell| \geq |r|$ for all
rules $(\ell,r) \in S$.

Now, consider an input $w\in \GG^*$ of length
$n=|w|$. For a moment, let a \emph{configuration} be  a pair $(u,v)$ such
that (i) $w \DAS*S uv$ and (ii) $u$ is irreducible. The goal
is to transform the initial configuration  $(1,w)$
in $\Oh(n)$ steps into some final configuration  $(\hat w, 1)$.

Say, we are in the configuration $(u,v)$. The goal is achieved
if $v= 1$. So assume that $v = av'$ where $a$ is a letter.
If $ua$ is irreducible, then we
replace $(u,av')$ by $(ua,v')$, and $(ua,v')$ is the next configuration. If however $ua$ is reducible, then
we can write $ua= u'\ell$ for some $(\ell,r) \in S$; and $u'$ is irreducible.
So, we
replace $(u,av')$ by $(u',rv')$, and $(u',rv')$ is the next configuration. The algorithm is obviously correct.
Defining the \emph{weight} $\gamma$ of configurations by
$\gamma(u,v)= (1-\eps) |u| + |v|$ we see that $\gamma$ reduces
from one configuration to the next by at least $\eps$. Hence
we have termination in linear time.
\end{proof}

In fact, length reducing rewriting systems arise naturally in
the class of small cancellation groups, and more generally
hyperbolic groups, which we might regard as a paradigm for
groups with easily solvable word problem.  To be precise:
a group $G$ is hyperbolic if and only if there is a finite generating set
$\GG$ for $G$ and a finite length-reducing system $S \subseteq \GG \times \GG$ (so
$G = \GG/S$) such that a word $w$ represents the trivial element of $G$ if and  only if $w$ can be $S$-reduced to the empty word, see \cite{ABC}.  In other words, a group is hyberbolic, \IFF{} there exists a finite length-reducing system which is confluent on the empty word.

\begin{definition}
A length-reducing string rewriting system which is
confluent on the empty word
is called a {\em Dehn} system.
\end{definition}
 If a group is defined by a finite length-reducing Dehn rewriting system then the rewriting algorithm is known in group theory as the  {\em Dehn} algorithm. More general definitions of Dehn algorithms, to rewriting systems
over a 
larger alphabet than the generators of the group, have been studied by Goodman and Shapiro 
\cite{GoodShap08}
and Kambites and Otto \cite{KambitesOtto08}. In particular in \cite{GoodShap08} it is shown that such generalised 
Dehn algorithms solve the word problem in finitely generated nilpotent groups and many
relatively hyperbolic groups.  

It is known that, given a finite presentation of a hyperbolic group $G$,
one can produce a finite Dehn presentation of $G$ by adding,
to a given presentation, all new relators of $G$  up to some length (which depends on the hyperbolicity constant of $G$).
However, this algorithm is very inefficient and the following questions
remain.

\begin{problem}
Is there a Knuth-Bendix type completion process that, given a finite presentation of  a hyperbolic group $G$,  finds a finite Dehn presentation of $G$.
\end{problem}

\begin{problem}
Is there an algorithm that, given a finite presentation of a hyperbolic
group,
determines whether or not  this presentation is Dehn.
\end{problem}

Notice that some partial answers to this question are known.  Namely, in \cite{Arzhantseva} Arzhantseva has shown that there is an algorithm that,
given a finite presentation of a hyperbolic group and  $\alpha \in [3/4,1)$,  detects whether or not
this presentation is an $\alpha$-Dehn presentation.
Here a presentation $\langle X \mid R\rangle$ of a group $G$ is called
an $\alpha$-Dehn presentation  if any non-empty freely reduced word $w \in (X \cup X^{-1})^*$ representing the
identity in $G$ contains as a factor a word $u$ which is also a factor of a cyclic
shift of some $r \in R^{\pm 1}$ with $|u| > \alpha|r|$.

\subsection{Infinite length-reducing systems}

Let us discuss some algorithmic aspects of rewriting with infinite length-reducing systems.

\begin{proposition}
Let $ S \subseteq \Gamma^* \times \Gamma^*$ be an infinite  recursive string rewriting system. Then the following hold.
 \begin{enumerate}
 \item [1)] If $S$ is length-reducing then an irreducible descendant of a given word can be computed. 

 \item [2)] If $S$ is Dehn and $M_S$ is a group, then the word problem in  $M_S$ is decidable.
 \end{enumerate} 
 \end{proposition}
\begin{proof}
The system $S$ is effective since it is recursive and length-reducing (see remark before Proposition \ref{pr:eff}). Now the argument in the proof of Proposition \ref{pr:eff} shows that for a given $w$ one can effectively find an $S$-irreducible of $w$, so  1) and 2) follow. 
\end{proof}

In the case of length-reducing systems one can try to estimate the time complexity of the algorithms involved. To this end we need the following definition.  Let $S$ be an effective non-length increasing rewriting system and
$$
  \ell_0 \to r_0, \ell_1 \to r_1, \ldots , \ell_i \to r_i, \ldots 
  $$
 its standard enumeration  (see Section \ref{subsec:computing}). If  there an algorithm $\mathcal{A}$ and a polynomial $p(n)$ such that for every $n \in \mathbb{N}$ the algorithm $\mathcal{A}$ writes down the initial part of the standard enumeration   of $S$   with $|\ell_i| \leq n$ in time $p(n)$ then the system $S$ is called  {\em enumerable in time} $p(n)$ or {\em Ptime enumerable}. In particular, we say that $S$ is linear (quadratic) time enumerable if the polynomial $p(n)$ is linear (quadratic).

 \begin{proposition}\label{prop:ptime}
Let $ S \subseteq \Gamma^* \times \Gamma^*$ be an infinite non-length increasing string rewriting system, which is  enumerable in time $p(n)$. Then the following hold.
 \begin{enumerate}
 \item [1)] If $S$ is length-reducing then an irreducible descendant of a given word $w$ can be computed 
 in polynomial time. 

 \item [2)] If $S$ is Dehn and $M_S$ is a group, then the word problem in  $M_S$ is decidable in time 
 in polynomial time.  
 \end{enumerate}
 \end{proposition}
\begin{proof} Given a word $w$ one can list in time $p(|w|)$ all the rules $\ell \to r$ of the standard enumeration of $S$ with  $|\ell| \leq |w|$.  Now in time $\Oh(p(|w|)|w|^2)$ one can check whether one of the listed rules can be applied to $w$ or not.  This proves 1) and 2). 
 \end{proof}

\subsection{Weight-reducing systems}
Many results above can be generalised to weight-reducing systems.
A \emph{weight} $\gamma$ assigns to each generator $a$ a positive
integer  $\gamma(a)$ with the obvious extensions to 
words by $\gamma(a_1 \cdots a_n) = \sum_{i=1}^n \gamma(a_i)$. 
A system is
called {\em weight-reducing}, if for every rule $\ell \to r$
one has $\gamma(\ell) > \gamma(r)$. The following statements in this paragraph  are taken {}from  \cite{die87tr}. It is decidable
whether a finite system is weight-reducing by linear integer programming. The reason to consider weight-reducing systems 
is that there are monoids like $\oneset{a,b,c}^* / ab = c^2$
having an obvious finite convergent weight-reducing presentation, 
but where no finite convergent length-reducing presentation exists.
 
For groups the situation is unclear. Actually, the
following conjecture has been stated. 
\begin{conjecture}Let $G$ be a finitely generated group. 
Then the following assertions are equivalent: 
\begin{enumerate}[1)]
\item  $G$ is a \emph{plain} group, i.e., $G$ is a free 
product of free and finite groups. 
\item $G$ has a finite convergent length-reducing presentation. 
\item $G$ has a finite convergent weight-reducing presentation.
\end{enumerate}
\end{conjecture}
The implications $1) \implies 2) \implies 3)$ are trivial, 
and 
$2) \implies 1)$ is known as the \emph{Gilman conjecture} and was stated first in 
\cite{Gilman84}. 

It is clear that the conjugacy problem can be decided in plain groups and this
holds for   groups $G$ having a finite convergent weight-reducing presentation, too. 
In fact, for $s,t\in G$ the set $R_{s,t}= \set{g\in G}{gs g^{-1}= t}$ is an effectively computable  rational 
subset of $G$. 

\section{Preperfect systems}\label{sec:preperfect}

\subsection{General results}

In this section we discuss preperfect rewriting systems, which play an important part  in solving the  word problem and finding geodesics (shortest  representatives in the equivalence classes) in groups.

\begin{definition}\label{def:thue}
A \emph{Thue system} is a rewriting system
$S \subseteq \Gamma^* \times \Gamma^*$ such that the following conditions hold:
\begin{enumerate}[i)]
\item If $\ell \longrightarrow r \in S$ then
$\abs{\ell} \geq \abs{r}$.
\item If $\ell \longrightarrow r \in S$ with $\abs{\ell} =
\abs{r}$,
then  $r \longrightarrow \ell \in S$, too.
\end{enumerate}
\end{definition}

%

To every rewriting system is associated an equivalent Thue system.
In order to specify a Thue system which is equivalent to a rewriting system
$S$ one can do the following:   symmetrize $S$ by adding all the rules $r \longrightarrow \ell$
whenever  $\ell \longrightarrow r \in S$, then throw out
all the length increasing rules. The new system, denoted $T(S)$  is called the
{\em Thue resolution}
of  $S$. It follows that every monoid  has a Thue presentation.

\begin{definition}\label{def:preperfect}
A confluent Thue system is called \emph{preperfect}.
\end{definition}

The main interest in preperfect systems in algebra comes from the following known (and easy) complexity result: for which we require the following definition.
\begin{definition} 
A word $w \in \Gamma^*$ is termed  $S$-{\em geodesic},\label{page:geodesic}
with respect to a string rewriting system
$S$, if it has minimal length in its $\DAS{*}{S}$-equivalence class
(and simply {\em
geodesic} where no ambiguity arises).
\end{definition} 
Clearly, $S$-geodesic words 
are precisely the geodesic words in the monoid $\GG/S$ relative to the generating set $\Gamma$, i.e.,  they have minimal length among all the words in $\Gamma^*$ that represent the same element in $\GG/S$. Sometimes, we say that a word $w \in \Gamma^*$ is a geodesic of a word $u \in \Gamma^*$,  if $w$ is $S$-geodesic and $\DAS * S$-equivalent to $u$.

\begin{proposition} \label{prop:preperfect-complexity} If a rewriting system $S$ is finite and preperfect, then
one can decide the word problem in the monoid defined by $S$
in polynomial space, and hence in exponential time. Moreover, along the way one can find an $S$-geodesic of a given word $w$, as well as, all $S$-geodesics of $w$.
\end{proposition}

A locally confluent  (strictly) length-reducing
system is convergent, hence, from the above, preperfect.
However the Thue resolution of an arbitrary finite convergent rewriting
system may fail to be terminating or confluent as simple examples show.
(Let $\Gamma=\{a,b,c,d,u,v\}$ and $S$ be the system with rules
$ab\ras{}{} u$, $bc\ras{}{} v$, $uc\ras{}{}d^3$ and $av\ras{}{}d^3$. Then
$T(S)$ is not confluent. The system with one rule $a\ras{}{}b$ has non-terminating Thue resolution.) It is also easy to see that $T(S)$ may be preperfect
when $S$ is not confluent.
On the other hand,  if a confluent system $S$ has no length-increasing rules, then the Thue resolution  can be constructed by
symmetrizing $S$ relative to all length preserving rules in $S$ (by adding the rule $r \longrightarrow \ell$ for each length preserving rule $\ell \longrightarrow r \in S$)  and a straightforward argument shows that in
this case $T(S)$ is confluent, so  preperfect.

 \begin{lemma}
If $S$ is a confluent rewriting system with
no length-increasing rules then the Thue resolution $T(S)$ is preperfect.
\end{lemma}

For a system $S$  (for example a Thue system) where
all rules $\ell \to r \in S$ are either length-reducing $|\ell| > |r|$, or length-preserving $|\ell| = |r|$, it is convenient to split  $S$ into a length reducing part $S_R$
and a length preserving part $S_{P}$, so $S = S_R \cup S_P$. 
If $S$ is a Thue system then all $S$-geodesic words that lie in the same equivalence class have the same length and any two of them are $S_P$-equivalent (can be transformed one into other by a sequence of rules from $S_P$).  Therefore, the word length in $\Gamma^*$ induces a well-defined length on the factor-monoid $M=\Gamma^*/S_{P}$ (application of relations from $S_P$ does not change the length). Hence, one can view
$S_R$ as a length reducing rewriting system over the
monoid $M=\Gamma^*/S_{P}$, in which case we assume that  $S_R \subseteq M \times M$.
Note that if $S_P$ is finite
and if there is an effective way to
perform reduction steps with $S_R$, then
the word problem in $M$ is  decidable.

Decidability of the word problem in $M=\Gamma^*/S_{P}$ allows one to  test whether a given rule from  $S_R$ is applicable to an element
of $M$. Since $S_R \subseteq M \times M$ is terminating it suffices to show
local confluence to ensure convergence. This may tempt one to introduce   an analogue of the Knuth-Bendix completion. However,
in general an infinite number of \emph{critical pairs} may appear in the
Knuth-Bendix process, and one needs to be able to recognise when
the  current system becomes preperfect. Unfortunately,
this is algorithmically undecidable.
More precisely, the following result holds.

\begin{theorem}[\cite{NarendranMcNaughton84}] \label{thm:prepundec}
The  problem of  verifying whether a finite Thue system is
preperfect or not  is undecidable.
\end{theorem}

In fact  in \cite{no88} this problem is shown to be undecidable even in the
case of a Thue system, whose length-preserving part $S_P$  consists only of   a single rewriting rule of the
form $ab\das{}{}ba$. On the other hand,
under some additional assumptions such a procedure can yield
useful results (\cite{die89tcs,die90}) - good examples in our context
are graph groups, c.f. Section~\ref{sec:gg}.

In the final part of this  section we discuss some complexity  issues in computing with  preperfect systems.
By Proposition \ref{prop:preperfect-complexity} finite preperfect systems allow one to solve the word problem and find geodesics in at most exponential time.

\begin{proposition}
Let $S$ be an infinite  preperfect rewriting system. Then:
 \begin{enumerate}
 \item [1)]  if $S$ is recursive then the word problem in the monoid $M_S$ defined by $S$ is decidable;

 \item [2)] if $S$ is Ptime enumerable then one can solve  the word problem in $M_S$ and find a geodesic of a given word in  exponential time.
 \end{enumerate} 
 \end{proposition}
\begin{proof}
Since preperfect rewriting systems are non-length-increasing it follows that recursive preperfect systems are effective. (see the remark before Proposition  \ref{pr:eff}). Therefore, given a word $w$ one can effectively list all the rules in $S$ with the left-hand sides  of length at most $|w|$. Denote this subsystem of $S$ by $S_w$. Now rewriting  $w$ using $S$ is exactly the same as using $S_w$, so 1) and 2) follow from the argument in the proof of Proposition \ref{prop:preperfect-complexity} for finite preperfect systems.

    \end{proof}

\section{Geodesically perfect rewriting systems}

In this section we consider a subclass  of Thue systems
which are designed to deal with geodesics in groups or monoids. In particular, we study confluent geodesic systems, which form  a subclass of  preperfect string rewriting systems, and which behave better
in many ways than general preperfect systems.
We call these systems \emph{geodesically perfect}, as this indicates
their essential properties and fits with the terminology of preperfect
systems.
However as discussed in Section \ref{sec:intro} they
are also known in the literature as
almost confluent or quasi-perfect.  The motivation for the  study
of geodesically perfect systems in group theory  comes mainly from attempts to
solve the, algorithmically difficult,
{\em geodesics problem}: that is, given a finite presentation of
a group $G$ and a word $w$ in the generators, find a word of
minimal length representing $w$ as an element of  $G$.

\subsection{Geodesic systems}\label{sec:geodesic}
We consider first a somewhat
larger, less well-behaved, class of rewriting systems.
\begin{definition}
A string rewriting  system $S \subseteq \Gamma^* \times \Gamma^*$ is called {\em geodesic} if   $S$-geodesic
words 
are exactly those words to which no length reducing rule from $S$ can
be applied.
\end{definition}

Note that if  $S$ is a  geodesic rewriting system then its Thue resolution is
also geodesic, this allows us to assume, without loss of generality,  that geodesic systems are  Thue systems.

\begin{remark}  Dehn rewriting systems are not in general geodesic: they  need only rewrite words that represent the identity   to (empty) geodesics in $\Gamma^*/S$.
\end{remark}

A finite geodesic system gives a linear time algorithm to find a geodesic of a given word $u \in \Gamma^*$.
The following algebraic characterisation of finite geodesic systems
in groups is given in \cite{GilHHR07} (the definition of geodesic in \cite{GilHHR07} is slightly
more restrictive than ours, however this makes no difference to the result).

\begin{theorem}[\cite{GilHHR07}]\label{th:fin-geodesic}
A group $G$ is defined by a
finite geodesic system $S$ \IFF{} $G$ is a finitely generated virtually free
group.
\end{theorem}

From the result of Rimlinger quoted above finitely generated virtually
free groups  are precisely
the universal groups of  finite pregroups.
It follows that every finite length reducing
geodesic system can be transformed to the length reducing  part of the rewriting system
(see Section \ref{sec:rewrite-U(P)}) associated with a finite pregroup.

The following result follows from Proposition \ref{prop:ptime}. 
\begin{proposition}

  \label{prop:comp-geod}
  Let $S$ be a geodesic Ptime enumerable string rewriting system such that the monoid $M_S$ is a group.  Then the word problem in the group $M_S$  is decidable in polynomial time.
  \end{proposition}

Very little is known about geodesic systems which do  not present
groups.
In particular, it is not clear whether the word problem remains decidable:
that is, given $u,v\in \Gamma^*$ decide whether or not  $u \DAS*S v$.
\begin{problem}\label{nixWP}
Does there exist a finite geodesic system $S$ for  which  the word problem
is undecidable?
\end{problem}

The following result demonstrates
one of the principal  difficulties of working with geodesic
systems.
\begin{theorem}\label{geoundec}
It is undecidable whether 
a finite rewriting system is geodesic.
\end{theorem}

\begin{proof}
The proof is a modification of the proof by Narendran and Otto \cite{no88}
which showed undecidability of preperfectness in
presence of a single commutation rule.

We need some notation and we adhere as far as possible to
that of \cite{no88}. We shall define the computation of a Turing
machine by a set of rewriting rules. A configuration of the machine
is then a particular form of word over the tape alphabet, the states and
the end markers. In detail let $\SIG$ be a finite set, the
tape alphabet, let $\ov \SIG$  be a disjoint copy of $\SIG$, let
$Q$ be a finite set of states, and $\alpha$ and $\beta$ be
special symbols representing end markers. There are two marked states $q_0$ and $q_f$, the
initial and final states.
The computation of the machine can be described by  a finite
   set
  of rules which fall into the following categories,
  where we use the notation $p,q \in Q$, $p\neq q_f$, $a, a', b \in \SIG$:
\begin{enumerate}[1.)]
\item  $pa \ras{}{} \ov {a'} q$.\\ (Read $a$ in state $p$, write
$a'$, move one step to the right, switch to state $q$.)
\item  $\ov b pa \ras{}{}  q b a'$. (As above, but move one step to the left.)
\item $p\beta \ras{}{}  pa\beta$. (Create new space before
the right end marker.)
\end{enumerate}
These rewriting rules constitute the rewriting system associated to $M$.
We assume that the machine is deterministic, so there no overlapping rules.
A {\em configuration} of a (deterministic) Turing machine is then  a word
$\alpha \ov u q v\beta $ with  $\ov u \in \ov{\SIG}^*$, $v \in \SIG^+$, and
$q \in Q$.
The initial configuration on input
$x \in \SIG^*$ is the word
$\alpha q_0 x\beta$. We assume that the machine stops \IFF{}
it reaches  the state $q_f$.

Now let $M$ be a Turing machine
for  which
it is undecidable whether or not computation halts on
on input $x \in \SIG^*$.
Using this machine we
are going to construct, for each $x\in \SIG^*$,
a new length reducing rewriting system $S_x$,
which is geodesic
\IFF{} the machine $M$ does not stop on input $x$.

The alphabet $\Gamma$ of each such system is to consist of
the symbols of $\SIG\cup \ov{\SIG}\cup Q\cup\{\alpha,\beta\}$
and  new additional symbols $d,e, \gamma, \delta, I,C$.
The system $S_x$ will consist of rules which simulate the computation
of $M$ on input $x$, with some additional control on the number of
steps of the computation carried out. Let $x\in \SIG^*$ .
To begin with, we introduce rules leading to two different
initial configurations. Let $m = \abs{x} + 5$. We define the two rules
$$\alpha q_0 x \beta\gamma \las{}{(1)} IC^m  \ras{}{(2)}\alpha q_0 x \beta \delta.$$

Next we introduce rules, involving  $d,e, \gamma$ and  $\delta$, to  control
the number of steps of the simulation.
Symbols $\gamma$ and $\delta$
convert $d$'s to $e$'s. The latter act as tokens to control the
number of steps performed by the simulation of $M$.
Both $\gamma$ and $\delta$ move right  consuming three $d$'s and producing
two $e$'s,
the difference being that between
$\gamma$ may move to the right arbitrarily far from $\beta$
whereas
$\delta$
is forced to remain very close to $\beta$. The effect on the
length of a word of each rule is the same.

Explicitly, we add new rules of the form:
$$ \gamma ddd \ras {}{}ee \gamma, \quad
\beta\delta ddd \ras {}{}\beta ee \delta.$$
Note that, using rule (1),
all words in $IC^m (ddd)^* $ now reduce as follows
$$IC^m d^{3n}  \RAS{}{(1)} \alpha q_0 x \beta\gamma d^{3n}
\RAS*{} \alpha q_0 x \beta  e^{2n}\gamma .$$
However, using  the rule $IC^m  \ras{}{(2)}\alpha q_0 x \beta\delta$
in the first step
we  can only do
$$IC^m d^{3n} \RAS{}{(2)} \alpha q_0 x \beta\delta d^{3n}
\RAS*{} \alpha q_0 x \beta  ee\delta d^{3n-3}$$
and then, for $n\geq 2$, we are stuck.

Now we bring $e$ into the game. The letter $e$ is used
to  \emph{enable} a 
computation step of $M$. It can move to the left until it
is at distance one to the right of
a state symbol.
The \emph{generic} rules for $e$ allow $e$ to move left and are as follows:
$$  abee\ras{}{} aeb, \quad aebee\ras{}{} aeeb \quad \text{for } a\in \SIG,
b \in \SIG \cup \{\beta\}.$$

Let us describe the effect of these rules on words
of the form
$$\alpha \ov u p a_1 \cdots a_k \beta \delta d ^{3n}$$
where $n$ is huge (and $k$ is viewed as constant $ k \geq 0$), $a_i \in \Sigma$, $\ov u \in {\ov \SIG}^*$.
The maximal possible reduction leads to a word of the form
$$\alpha \ov u p  a_1 ee \cdots a_k ee \beta ee \delta d ^{3n'}.$$
In this case, if  $n$ is large enough, then  $n'> 0$, no further reduction
is possible
and actually $n-n' \in \Oh(1)$.

At this point we introduce rules to simulate the computation
of the machine $M$.
There is one simulation rule
corresponding to each rule in the
rewriting system associated to $M$.
More precisely we introduce a rule
$uee\ras{}{} v$ for  each rewriting rule $u\ras{}{} v$ of $M$: so we
have simulation
rules of three types
(where again we use the notation $p,q \in Q$, $a, a', b \in \SIG$)
$$ paee \ras{}{} \ov {a'} q, \quad
\ov b paee\ras{}{}  q b a', \quad
p\beta ee \ras{}{}  pa\beta.$$
The system $S_x$ consists of the rules defined so far, which we
list in Figure \ref{fig:sysx}, so is length reducing.
\begin{figure}[bt]
\begin{enumerate}[I.)]
\item Initial rules:
$$\alpha q_0 x \beta\gamma \las{}{(1)} IC^m  \ras{}{(2)}\alpha q_0 x \beta \delta.$$
\item Step control rules, for $a\in \SIG$, $b\in \SIG\cup \{\beta\}$:
\begin{align*}
\gamma ddd & \ras {}{}ee \gamma,\\
\beta\delta ddd & \ras {}{}\beta ee \delta, \\
abee& \ras{}{} aeb, \\
aebee& \ras{}{} aeeb.
\end{align*}
\item Simulation rules, for $a,b\in \SIG$, $p\in Q\backslash\{q_f\}$:
\begin{align*}
paee & \ras{}{} \ov {a'} q, \\
\ov b paee& \ras{}{}  q b a', \\
p\beta ee & \ras{}{}  pa\beta.
\end{align*}
\end{enumerate}
\caption{The system $S_x$.}
   \label{fig:sysx}
 \end{figure}

Now assume that the machine $M$ halts on input $x$.
This implies that  only finitely many computation steps $t$ can be performed.
Again choose $n$ huge and view  $t$ and $\abs{x}$ as constants.
Consider a word of the form
$IC^md^{3n}.$
Starting a  reduction with the second rule we get stuck at an
irreducible word when the
simulation reaches state $q_f$:
$$IC^md^{3n} \ras{}{(2)}
\alpha q_0 x\beta \delta d ^{3n} \RAS{*}{}\alpha \ov u q_f a_1 ee \cdots a_k ee \beta ee \delta d ^{3n'}$$
at which point $n-n'\in \Oh(1)$.
The system $S_x$ cannot be geodesic because with the other initial rule we
can first move $\Gamma$ to the right of all the $d$'s thereby losing
$n$ letters immediately:
$$IC^md^{3n} \ras{}{(1)}
\alpha q_0 x\beta \gamma d^{3n} \RAS{*}{}\alpha q_0 x\beta d^{2n}\gamma $$
and then when
the simulation reaches the state $q_f$ the resulting irreducible word will
end $e^{2n'}\gamma$ instead of $\delta d^{3n'}$ (and otherwise will
be the same).

It remains to cover  the case when the machine does not
halt on input $x$. We shall  show that in this case
the system $S_x$ is geodesic. Note that, as $M$ never reaches
state $q_f$, for all $n>0$
$$\alpha q_0x\beta e^{2n}  \RAS*{S_x}
\alpha \ov{u} p y \beta,
$$
where $\ov{u}\in \ov{\SIG}^*$, $p\in Q$ and $y\in (\SIG \cup \{e\})^*$.
For technical reasons we define a sequence of
words $w_i$ for $i \geq 0$ as follows. We let $w_0 = q_0x$
and let $\alpha w_{i+1}\beta$ be defined to be the irreducible descendant of
$\alpha w_{i} \beta ee$.
The sequence of words $w_i$ is infinite because the machine does not stop
on input $x$.
Thus
$$\forall i\geq 0: \; \alpha w_{i} \beta ee\RAS*{S_x}\alpha w_{i+1} \beta\in \IRR(S_x).$$

Now we add infinitely many rules to $S_x$ to form a new system
$T_x$ as follows:
$$\forall i\geq 0: \; \alpha w_{i} \beta  \delta \ras{}{}
\alpha w_{i} \beta  \gamma
$$

As the rules of $T_x$ are generated
by steps of the  Knuth-Bendix completion  procedure applied to $S_x$ the
congruences generated by $S_x$ and $T_x$ are the same.
To summarise, the system $T_x$ consists of the rules of Figure \ref{fig:sysx}
and those listed in Figure \ref{fig:syst}. Thus $T_x$ is
is terminating and local confluence can be checked directly.
\begin{figure}[bt]
\begin{enumerate}[I.)]
\setcounter{enumi}{3}
\item Completion rules:
$$ \forall i\geq 0: \; \alpha w_{i} \beta  \delta \ras{}{}
\alpha w_{i} \beta  \gamma
$$ 
\item[]where
$\alpha w_{i} \beta ee\RAS*{S_x}\alpha w_{i+1} \beta, \forall i\geq 0,\;
\textrm{ and } w_0=q_0x.$
\end{enumerate}
\caption{The additional rules of system $T_x$.}
   \label{fig:syst}
 \end{figure}

Each word $w\in \Gamma^*$
has a unique factorisation where we choose $k$ and all $n_j$ to be maximal:
$$ u_0 (IC^m d^{n_1}) u_1 \cdots (IC^m d ^{n_k}) u_k.$$
The benefit of the system
$T_x$ is that it provides us with canonical geodesics.
A geodesic of $w$ is given by:
$$ \widehat{u_0} (\alpha w_{i_1} \beta \gamma d^{m_1})
\widehat{ u_1} \cdots (\alpha w_{i_k} \beta\gamma d^{m_k}) \widehat{u_k},$$
where $m_j=n_j \bmod 3$.
The crucial observation is that
allowing only rules from $S_x$ we achieve exactly the same
form with the exception that some $\gamma$'s are still $\delta$'s.  Thus, the system is geodesic.
 \end{proof}

\subsection{Geodesically perfect systems}\label{sec:geodesicperf}

\begin{definition}
A string rewriting system  $S \subseteq \Gamma^* \times \Gamma^*$ is called
{\em geodesically perfect},  if
 \begin{enumerate}  [i)]
 \item $S$ is geodesic, and
 \item if $u, v \in \Gamma^*$ are $S$-geodesics, then  $u \DAS * S v$ if and
only if $u \DAS * {S_{P}} v$, where  $S_{P}$ is the length-preserving part of
$S$.
     \end{enumerate}
 \end{definition}

Again, it follows directly that if $S$ is a geodesically perfect system
then so is its Thue resolution, so we can assume that geodesically perfect systems are Thue.
If $S$ is a geodesically perfect  Thue system then
we write it as  $S = S_R \cup S_{P}$ where
$S_R$ is its length reducing part and $S_{P}$ its  length preserving part.
It also follows from the definition that a geodesically perfect system is
confluent.

There is a simple procedure to describe geodesics of elements in  the monoid
$\Gamma^*/S$ defined by a \gp{} Thue
system $S$. Namely, the geodesics of a given
word $w \in \Gamma^*$ are the $S_R$-reduced forms of $w$ and any two such
geodesics can be obtained from one another by applying finitely many  rules
from $S_{P}$. Moreover it is shown in \cite{bo93springer} 
that  the word problem for monoids defined by finite geodesically perfect
rewriting systems is PSPACE complete. 

The following result relates geodesically perfect to preperfect Thue systems.
\begin{proposition} \label{pr:geodesically perfect}
Let $S \subseteq \Gamma^* \times \Gamma^*$ be a Thue system. Then
\begin{enumerate} [1)]
 \item if $S$ is geodesically perfect then it is preperfect and
 \item if $S$ is preperfect and  geodesic then it is geodesically perfect.
\end{enumerate}
\end{proposition}
\begin{proof}
1) follows from the observation that geodesically perfect implies
confluent. To see 2) observe that $S$ is confluent, hence Church-Rosser.
Therefore, if  $u, v$ are two geodesics with $u \DAS * S v$ then $u \RAS * S
w$ and $w \LAS * S v$ for some $w \in \Gamma^*$. Since $u, v $ are
$S$-geodesics the only rules that could be applied in $u \RAS * S w$ and
$w \LAS * S v$ are length preserving, hence $u \DAS * {S_P} v$, as required.
\end{proof}

In Section \ref{sec:rewrite-U(P)} we will describe a general tool to
construct geodesically perfect systems defining groups: based on the
fact that rewriting systems  associated with pregroups are  always
geodesically perfect.

In Corollary~\ref{virtfreegp} we prove that groups defined by finite geodesic systems are  exactly the groups defined by finite geodesically perfect  systems.

Obviously, every geodesic rewriting system $S$ contains
the length-reducing  part $T_R$ of some (infinite) geodesically perfect Thue
system $T$ defining the same monoid.
Indeed, one can obtain $T$ by first constructing the Thue
resolution $T^\prime$ of $S$ and then adding  length-preserving rules to
$T^\prime$ to make it confluent. But it is not true that every finite geodesic rewriting system $S$ is the length-reducing  part of a finite geodesically perfect system  defining the same monoid.  To see this consider the following example.

\begin{example}\label{geoper}
The following system is geodesic, and it is not the  length-reducing part of any finite geodesically perfect system
defining the same quotient monoid.
$$ add \ras{}{} ab,\quad  add \ras{}{} ac, \quad
  bdd \ras{}{} eb,\quad  cdd \ras{}{} ec.$$
Indeed, let $S$  be the system above, and let
$T = S \cup \oneset{b \longleftrightarrow c}$.
The new system $T$ is \gp{} by
Proposition~\ref{prop:nivat}. But $T$-geodesics are computed by
using rules {}from $S$. As $\DAS*S \subseteq \DAS*T$ we see  that  $S$ is a geodesic system.

Let us show that $S$ is not the  length-reducing part of any equivalent,
finite,
geodesically perfect system.
For a contradiction, assume that
a finite set $T$ of non-trivial symmetric rules can be added to
$S$ such that $S\cup T$ becomes \gp and is equivalent to $S$.
Assume $T$  involves a new letter,
say $f$. Then $f$ is equal
to some word $u_f$ over $\{a,b,c,d,e\}$
which is irreducible with respect to $S$. If  $u_f$ is empty, then we do not
need $f$, hence $u_f$ is nonempty and
we have $u_f\RAS*T f$. The rules of $T$
are symmetric (hence length preserving),
so   $f$  is
accompanied by a rule, say $f \longleftrightarrow a$, and $f$ is redundant.
So,  actually
we may assume $ T \subseteq \{a,b,c,d,e\}^* \times \{a,b,c,d,e\}^*$.
Clearly, $ae^n b \LAS*{S}ad^{2n+2}\RAS*{S}ae^n c$,  hence
$ae^n b \DAS*{S} ae^n c$ and  $ae^n b$, $ae^n c$ are in the same class and
are $S$-reduced.  Because $T$ is finite,  some left hand side of $T$ must
contain
a word in $u \in ae^* \cup e^* \cup e^*b \cup e^*c$.
But all these words $u$  are $S$-reduced, hence geodesic.
Moreover, for any such $u$ there is no other word $v$ in the same class
as $u$ and of the same length. So, for large enough $n$ the rules
of $T$ cannot be applied to either $ae^n b$ or $ae^n c$. As
$T$ is the length preserving part of the supposedly \gp{} system $S\cup T$, this
is the required contradiction.
 \end{example} 

\begin{remark}\label{rem:gpex}
Let $M = \{a,b,c,d,e\}^*/S$ be the quotient monoid as in Example~\ref{geoper}.
The proof above can be modified in order to show that
actually there is no finite system $ T \subseteq \{a,b,c,d,e\}^* \times \{a,b,c,d,e\}^*$
which is \gp{} 
and which defines $M$. However, if we use an additional letter $f$,
then the following
system defines $M$, too.
$$ dd \ras{}{} f,\quad
af \das{}{} ab,\quad  af \das{}{} ac, \quad
  bf \das{}{} eb,\quad  cf \das{}{} ec.$$
  The system is \gp{}, by Proposition~\ref{prop:nivat} again.
\end{remark}

We note that Example~\ref{geoper} illustrates a general fact: namely
that if $S$ is a rewriting system and there exists a set $T$ of symmetric
rules such that $S\cup T$ is geodesically perfect (but not necessarily
equivalent to $S$) then $S$ itself is geodesic. Since geodesic systems
are undecidable whereas geodesically perfect systems are decidable this
could prove to be a useful test for a geodesic system.

\section{Knuth-Bendix completion for geodesically perfect systems}\label{sec:KB}

A classical result of Nivat and Benois 
(stated in Proposition~\ref{prop:nivat})
shows that it is decidable whether a finite Thue system is geodesically perfect.
In order to explain the criterion we need the notion of
\emph{critical pair}. All rewriting systems $S$ in this subsection
are viewed as Thue systems and split into a length reducing part $S_R$
and a length preserving part of symmetric rules $S_P$.
By definition, a critical pair is a pair $(x,y)$ arising from the situation
$$x \LAS{}{(\ell_1,r_1)} z \RAS{}{(\ell_2,r_2)} y$$
subject to  the following conditions:
\begin{enumerate}[1.]
\item    $(\ell_1,r_1)\in S_R$ is length reducing, but
$(\ell_2,r_2)\in S$ can be any rule.
\item $z = \ell_i u_i = u_j \ell_j$ with $|u_i| < |\ell_j|$ and $i,j \in \{1,2\}$ such that $i = j$ implies $u_i= u_j = 1$.
\end{enumerate}

\begin{proposition}[\cite{nibe71}]\label{prop:nivat}
A finite Thue system $S$ is \gp{} \IFF{} for all critical pairs
$(x,y)$ there are words $ x'$ and $ y'$ such that
with length reducing reductions we have:
$$ x' \LAS*{S_R} x, \quad y \RAS*{S_R} y',$$
and with length preserving reductions we have:
$$  x' \DAS*{S_P} y'$$
\end{proposition}

\begin{proof}The proof is not very difficult and can be found,
for example, in the book \cite[Thm.~3.6.4]{bo93springer}.
\end{proof}

\begin{remark}\label{rem:nibenote}
Note that the words $ x'$ and $ y'$ in Proposition~\ref{prop:nivat} need not 
be irreducible w.r.t. the length reducing subsystem $S_R$.
This fact is actually used in the proof
of Proposition~\ref{prop:kbcgp}.
\end{remark}

This criterion leads to the following version of the Knuth-Bendix procedure.
Consider a  finite Thue system $S_0$. We shall construct a series
of Thue systems $S_0 \subseteq S_1 \subseteq S_2 \subseteq \cdots $
such that the union over all $S_i$ is \gp{} and we have
$S_i = S_{i+1}$ (i.e., the completion procedure stops) \IFF there exists
a finite Thue system $T$ which is \gp{} and equivalent to $S_0$, that is ${}\DAS*{S_0}{} = {}\DAS*{T}{}$.
We divide the procedure into phases. We assume that in phase $i$
a Thue system $S_i = S_{R} \cup S_{P}$ has been defined such
that $S_{R}$ contains the length reducing rules, $S_{P}$
contains the length preserving rules, and ${}\DAS*{S_0}{} = {}\DAS*{S_i}{}$.

We begin phase $i+1$ by computing a list of all critical pairs
of the system $S_i$ (which were not already considered in phases $1$ to $i$).
For each such pair $(x,y)$ choose words
$\widehat x, \widehat y$, irreducible with respect to
the subsystem $S_R$, such that
$$ \widehat x \LAS*{S_R} x, \quad y \RAS*{S_R} \widehat y.$$
Define new rules as follows.
\begin{itemize}
\item If $|\widehat x|>|\widehat y|$ then add the rule
$\widehat x \ras{}{} \widehat y$ to $S_R$.
\item If $|\widehat y|>|\widehat x|$ then add the rule
$\widehat y \ras{}{} \widehat x$ to $S_R$.
\item If $|\widehat y|=|\widehat x|$ then  test whether or not
$$\widehat x \DAS{*}{S_P} \widehat y.$$
If the answer is negative then  add the symmetric  rule
$\widehat x \das{}{} \widehat y$ to $S_P$.
\end{itemize}
The system $S_{i+1}$ is defined to be  $S_{i}$ together with all new
rules which have been added to resolve all critical pairs of $S_i$.
On a formal level we define $S_i$ for all $i\geq 0$, but, of course,
the procedure stops as soon as $S_i = S_{i+1}$, i.e.,  no new
rules are needed to resolve critical pairs of $S_i$. Thus, if
it stops with $S_i = S_{i+1}$ then $S_i$ is a finite \gp{} Thue system,
which is equivalent to $S_0$ (and we have $S_i = S_{j}$ for all $i \leq j$).
However, what we really wish is stated in the
following proposition.

\begin{proposition}\label{prop:kbcgp}Let
$S_0 = S_{R} \cup S_{P}$ be a finite Thue system  with length reducing rules
$S_{R}$ and length preserving rules $S_{P}$. Let
$$S_0 \subseteq S_1 \subseteq \cdots S_i \subseteq \cdots $$
be the sequence of Thue systems which are computed by the
Knuth-Bendix completion as described above. Let $\widetilde S = \bigcup_{i\geq 0} S_i$. Then the system $\widetilde S$ is \gp{} and
we have ${}\DAS*{S_0}{} = {}\DAS*{S_i}{} = {}\DAS*{\widetilde s}{}$
for all $i \geq 0$. Moreover the following statements are equivalent.
\begin{enumerate}[1.)]
\item We have $S_i = S_{i+1}$ for some $i\geq 0$.
\item The Thue system $\widetilde S$ is finite and \gp.
\item There exists some finite \gp{} Thue system $T$ such that
${}\DAS*{S_0}{} = {}\DAS*{T}{}$.
\end{enumerate}
\end{proposition}

\begin{proof}
If $S_i = S_{i+1}$ for some $i\geq 0$, then clearly
$\widetilde S$ is finite. It is \gp{} by the criterion of Nivat and Benois, c.f.
Proposition~\ref{prop:nivat}. So, assume there exists
some  finite \gp{} Thue system $T$ with
${}\DAS*{S_0}{} = {}\DAS*{T}{}$. We have to show that the procedure stops.
We let $m$ be large enough that $m \geq  \max\set {|\ell|}{(\ell, r )  \in T}$.
Next we consider $i$ large enough  that $S_i$ contains all rules {}from
$\widetilde S$ where the left hand side has length of at most $m$.
Clearly, an index $i\in \NN$ with this property exists.
We will show that $S_i$ is \gp, and  Proposition~\ref{prop:nivat} immediately implies
$S_i = S_{i+1}$. For technical reasons, in a first step we
we remove {}from $T$ all  length preserving rules
$(\ell,r) \in T$
where we can apply to $\ell$ a length reducing rule of $T$.  It is clear that new and smaller system $T'$ is still \gp{} and
${}\DAS*{S_0}{} = {}\DAS*{T'}{}$; so we replace $T$ with
$T'$.
Since now $(\ell,r) \in T$
with $|\ell|= |r|$ implies that $\ell$ and $r$ are geodesics and since
$\widetilde S$ is \gp and $m$ is large enough,
we see that ${}\DAS*{T_P}{} \subseteq{}\DAS*{(S_i)_P}{}$.

Next consider some word $\widehat x$, which is irreducible with respect to the
length reducing rules in $S_i$. The claim is that $\widehat x$
is a geodesic. Indeed assume the contrary. Then a length reducing rule
$(\ell,r) \in T$ can be applied to  $\widehat x$. Since $\ell$ is not geodesic,
there is  a length reducing rule in $\widetilde S$ which can be applied to
$\ell$, but due to the definition of $m$ this rule is in $S_i$, too.
Thus, we have a contradiction and so $S_i$ is geodesic.
Now suppose that  $\widehat x$ and $\widehat y$  are geodesic
and  that $\widehat x \DAS{*}{S_i} \widehat y$.
Then $\widehat x \DAS*{T}\widehat y$ so $\widehat x \DAS*{T_P}\widehat y$.
As ${}\DAS*{T_P}{} \subseteq{}\DAS*{(S_i)_P}{}$ this implies
$\widehat x \DAS*{(S_i)_P}\widehat y$ so
$S_i$ is \gp.
\end{proof}


A finite  geodesic (or geodesically perfect)  rewriting system
$S \subseteq \Gamma^* \times \Gamma^*$ allows one to find $S$-geodesics
in linear time.
In particular, if the monoid $M = \Gamma^*/S$  defined by $S$  is a group
one can  solve the word problem in $M$ in linear time. However, in general, there seems to be no linear time reduction {}from
the word problem in a monoid  $M$  to the geodesic problem.

\section{Examples of preperfect systems in groups}
\label{sec:examples}

\subsection{Graph groups}\label{sec:gg}

Let $\Delta = (\Sigma, E) $ be an undirected graph. The
\emph{graph group} (or a \emph{right angled Artin group}, or a \emph{partially commutative group}) defined by $\Delta $ is the group $G(\Delta)$ given by the presentation
\begin{align*}
G(\Delta) = F(\Sigma)/\set{ab=ba}{(a,b)\in E},
\end{align*}
where $F(\Sigma)$ is the free group with basis $\Sigma$.
The group  $G(\Delta)$ has a monoid presentation given by  a preperfect
rewriting system $S_\Delta$. Indeed,  let $\Gamma= \Sigma \cup \ov{\Sigma}$
where $\ov{\Sigma}$ is a disjoint copy of ${\Sigma}$.
The rules of $S_\Delta$ are:
\[
\begin{array}{rcll}
a \ov{a} & \longrightarrow  &1 &  \; \\
ab& \longleftrightarrow  &ba&  \; \mbox{if } \oneset{(a,b), (\ov{a},b),
(a,\ov{b}), (\ov{a},\ov{b}) } \cap E \neq \emptyset
\end{array}
\]
where $a,b\in \Gamma$ and $\ov{\ov{a}} = a$ for all $a\in \Gamma$.

If the graph $\Delta $ is finite the system $S_\Delta $ provides us
with a decision algorithm for solving the word problem in $G(\Delta )$, though not the fastest one (WP in  graph groups can be solved in linear time, see \cite{wra88,die90}). However, the system $S_\Delta $  is very intuitive and simple, and  it gives the geodesics in $G(\Delta)$, which are precisely the words whose length cannot be reduced by $S_\Delta$.

Although it is preperfect the
system $S_\Delta$ is not geodesically perfect. However every graph
group may be constructed by a sequence of HNN-extensions and free
products with amalgamation, starting with infinite cyclic groups, and
so, from the results of Section \ref{sec:pregroup} below it follows
that these groups may be defined by (infinite) geodesically perfect systems.
Moreover  finite convergent rewriting systems for these groups have
been found by Hermiller and Meier \cite{HermillerMeier} (see also
\cite{BokutShiao,EsypKR,vanWyk}).

\subsection{Coxeter groups}

Let $D_3 = \{a,b\}^*/\{a^2= 1, b^2 = 1,(ab)^3 = 1\}$ be a dihedral group. Define a preperfect system $S$ by the following rules
\[
\begin{array}{lcl}
aa  & \longrightarrow  &1,   \\
b b & \longrightarrow  &1,  \\
aba& \longleftrightarrow  &bab.
\end{array}
\]

More generally,  a \emph{Coxeter group} on $n$ generators
$
a_1, \ldots, a_n
$ is given by a symmetric $n\times n$ matrix
$(m_{ij}$) with entries in $\NN$ and 1's on the diagonal. The defining
relations are given by:

$$ (a_ia_j)^{m_{ij}} = 1 \quad \text{ for all } 1\leq i,j \leq n.$$
Note that this implies $ a_i^{2} = 1$ since $m_{ii}= 1$; and if $m_{ij}= 0$,
then the equation  $ (a_ia_j)^{0} = 1 $ is trivial. (Therefore
it is also common to write $ (a_ia_j)^{\infty} = 1 $, because
$a_ia_j$ turns out to be an element of infinite order in this case.)

The word problem of Coxeter groups can be solved by the preperfect
\emph{Tits system} \cite{tits67}
(see also (\cite{Bourbaki02, Brown08, Davis08}) of rewriting rules:
\begin{align*}
a_i^2 & \longrightarrow  1,&&\text{for } 1\leq i \leq n,\\
(a_ia_j a_i a_j \cdots ) & \longleftrightarrow 
(a_ja_i a_j a_i \cdots ) &&\text{for } 1\leq i,j  \leq n \; \text{ and } \\
&&&\abs{(a_ia_j a_i a_j \cdots )} = \abs{(a_ja_i a_j a_i \cdots )} = m_{ij}.
\end{align*}

The classical proof that this system is preperfect relies on the fact that
Coxeter groups are linear \cite{bjofra05}. Of course this system
is not geodesically perfect. For virtually free Coxeter groups
Corollary \ref{virtfreegp} guarantees the existence of a finite
geodesically perfect rewriting system.
It is shown in \cite{GordonLongReid04} that every Coxeter group
is either virtually free or contains a surface group; but the
question of whether the latter can be
defined by a geodesically perfect system (necessarily infinite) remains open.

Convergent rewriting systems for Coxeter groups
have been constructed, using the
Knuth-Bendix procedure, by le Chenadec \cite{Chenadec86}, but in general these are
not finite. Finite convergent rewriting systems for certain classes
of Coxeter groups have been found by Hermiller
\cite{Hermiller94} (see also \cite{duCloux99,Borges04}).

\subsection{HNN-extensions}\label{sec:hnn}

Let $G$ be any group with isomorphic subgroups $A$ and $B$.
Let $\Phi: A \to B$ an isomorphism and
let $t$ be a fresh letter. By $\gen{G,t}$ we mean the free
product of $G$ with the free group $F(t)$ over $t$. The HNN-extension
of $G$ by $(A,B,\Phi)$ is the quotient group
\begin{align*}
\HNN(G;A,B,\Phi) = \gen{G,t}/\set{t^{-1}at=\Phi(a)}{a \in A}
\end{align*}
There is normal form theorem for elements in $\HNN(G;A,B,\Phi)$, which
implies that $G$ embeds into $\HNN(G;A,B,\Phi)$ and shows under which restrictions decidability
of the word problem for $G$ transfers to HNN-extensions.  Usually the
normal form theorem is shown by appeal to  a combination of arguments
of Higman, Neumann and Neumann and Britton,
see \cite[Chapter IV, Theorem 2.1]{LS}.

Another option is to define a convergent string rewriting system.
To see this, let
$\Gamma = \oneset{t,t^{-1}}\cup G\setminus \smallset{1}$ and
view $\Gamma$ as a possibly infinite alphabet.
We identify $1\in G$
with the empty word $1\in \Gamma^*$.
We choose transversals for cosets of $A$ and $B$. This means
we choose $X,Y \subseteq G$ such that there are unique decompositions
\begin{align*}
G= AX = BY
\end{align*}
We may assume that $1\in X\cap Y$.

The system $S \subseteq \Gamma^*
\times \Gamma^*$ is now defined by the following rules
with the convention that $[gh]$ denotes $gh\in G$ (as a single letter or the empty word).
$$
\begin{array}{ccll}
t^{-1} t & \longrightarrow  &1;  &
tt^{-1}   \longrightarrow  1;  \quad
gh \longrightarrow  [gh],\textrm{ for all } g,h\in G; \\
tg& \longrightarrow & aty,& \mbox{if } \;a\in A, \;a\neq 1, \;y \in Y,\; \Phi(a)y=g \mbox{ in } G;\\
t^{-1}g& \longrightarrow  & bt^{-1}x,& \mbox{if } \;b\in B,\; b\neq 1,\;
x \in X, \;  \Phi^{-1}(b)x=g \mbox{ in } G.
\end{array}
$$

\begin{proposition}\label{prop:hnnviaconvergence}
The system $S$ above is convergent and defines the HNN-extension of $G$ by
by $(A,B,\Phi)$.
Every irreducible normal form admits
a unique decomposition as
$$ 
g= g_0 t^{\varepsilon_1} g_1 \cdots t^{\varepsilon_n} g_n
$$ 
with $n$ minimal such that $n\geq 0$, $g_0 \in G\setminus \smallset{1}$,
and either
$\varepsilon_i = -1$ with $g_i \in X$ or $\varepsilon_i = 1$ with $g_i
\in Y$,
for all $1\leq i \leq n$.
\end{proposition}

\begin{proof}
Obviously, $\Gamma^*/S$ defines the HNN-extension of $G$ by
by $(A,B,\Phi)$.
Although the system has length-increasing rules it is not too
difficult to prove termination.  Local confluence is straightforward, so $S$ is indeed
convergent. Since all elements of $G$ are irreducible we see that
$G$ embeds into the HNN-extension. Moreover, it is also clear that we obtain the normal form as stated in the proposition.
\end{proof}

This convergent system also leads to the following well-known
classical fact.
\begin{corollary}\label{cor:hnnWP}
Assume that have the following properties:
$H$ is finitely generated and has a decidable word problem, membership problems for $A$ and $B$ are solvable, and the
isomorphism $\Phi: A \to B$ is effectively calculable.
Then the  HNN-extension of $G$ by
by $(A,B,\Phi)$ has a decidable word problem.
\end{corollary}

\begin{proof}
We may represent all group elements in $H$ by length-lexicographic
first elements (i.e., choose among all geodesics the  lexicographical first one).
The transversal
$X$ (resp. $Y$) may be chosen to consist of
the  length-lexicographic first  element
of each coset $Ag$ (resp. $Bg$), where $g$ runs over $G$.
Given $g$ we can compute
the representative of $Ag$ in $X$ (resp. $Bg$ in $Y$), because  membership is
decidable for $A$ and $B$. Now, given $b\in B$, the ability
to compute $\Phi$ allows us to find $a\in A$ with $\Phi(a)=b$.
Thus, all steps in computing normal forms are effective.
\end{proof}

It should be clear however that the purpose of the
system $S$ above is not to decide the word problem effectively; but rather
to facilitate straightforward proofs of other results, such as Britton's lemma.
Consider the following system $B$ of Britton reduction rules.
$$
\begin{array}{ccll}
t^{-1} t &\ras{}{}& 1; &
tt^{-1}  \ras{}{} 1; \quad
gh\ras{}{} f,\;   \mbox{ if } \; gh= f \; \mbox{ in } G;\\
t^{-1}a t & \ras{}{}& \Phi(a) &   \mbox{if } \;a\in A;\\
tb t^{-1}& \ras{}{} &\Phi^{-1}(b) &\mbox{if } \;b\in B\end{array}
$$

The system $B$ is length reducing, but not confluent. However,
$\RAS{}{B}\subseteq \RAS*H$, hence we can think of $B$ as
as subsystem of $H$.
Britton's lemma says that $B$ is confluent on all words
which represent $1$ in the HNN-extension.
Here is a proof using our system $S$.  Consider
any Britton reduced word $g$. It has the form
$ 
g= g_0 t^{\varepsilon_1} g_1 \cdots t^{\varepsilon_n} g_n.
$ 
Applying rules from $H$ does not destroy the property of
being Britton reduced and neither $t$ nor $t^{-1}$ can vanish.
Thus, if $g$ reduces to the empty word using $H$, then
$g$ is already the empty word.

Observe that $B$ is not a geodesic system, because
$a t \Phi(a)^{-1}$ is Britton reduced, but $a t \Phi(a)^{-1}= t$.
In Example \ref{ex:hnnpre} below we construct a geodesically perfect
rewriting system for an HNN-extension.

\subsection{Free products with amalgamation}
\label{sec:amal-prod-rewrite}

There is a natural  convergent (resp.{} \gp)
rewriting system which defines amalgamated products. Let
$A$ and $B$ be groups intersecting in a common subgroup $H$.
This time we  choose transversals for cosets of $H$ in $A$ and in $B$;
that is  $X\subseteq A$
and $Y\subseteq B$ with $1 \in X \cap Y$ such that there are unique decompositions
$A = HX$ and $B= HY$.
We let $\Gamma = (A \cup B) \setminus \{1\}$ and we identify
$1$ with the empty word in $\Gamma^*$.

We use the convention to write $[ab]$ for the product $ab$ whenever it is defined.
This means  $[ab]$ is viewed as a letter in $\GG$ or $[ab]=1 $ and it is defined if
either $a,b \in A$ or $a,b \in B$.

The system $S \subseteq \Gamma^2
\times (\{1\} \cup \Gamma \cup \Gamma^2)$ is now defined by the following rules:
$$
\begin{array}{rcll}
ab & \longrightarrow  &[ab] & \mbox{if $[ab]$ is defined,} \\
ab& \longrightarrow & [ah]y& \mbox{if } \;1 \neq a\in A,
\;h \in H, \;b \neq y \in  Y,\mbox{ and  } b = [hy]\\
ba& \longrightarrow & [bh]x& \mbox{if } \;1 \neq b\in B, \;h \in H,  \;a\neq x  \in X,\mbox{ and  } a = [hx]
\end{array}
$$
The system
defines the amalgamated product $G=A*_H B$.
It  is terminating by a length lexicographical ordering.
Local confluence follows by a direct inspection, whence
convergence.
Again we obtain the
normal form theorem (\textit{cf.} \cite[Corollary 4.4.1]{mks}):
every element $g$ of $G$ has a unique decomposition as
\[
g=[hg_0]g_1\cdots g_n,
\]
where $h\in H$, $g_i$ is a non-trivial element of $X\cup Y$ and $g_i$ and
$g_{i+1}$ do not lie in the same factor. However, in practice we may
not wish to compute transversals
explicitly. So let us apply only length reducing rules $ab  \longrightarrow  [ab]$ only until
we end up with a word
$ g = g_0\cdots g_n$, to which no length reducing
rule may be applied. Since we cannot apply length reducing rules to $g$ we obtain
that
$$  \forall 0 \leq i < n: \; g_i \in A \iff  g_{i+1} \in B\setminus H \; \wedge \;
g_i \in B \iff  g_{i+1} \in A\setminus H.$$
Further applications of the rules of $S$ preserve this property. Thus, $S$ is \gp, even if we use the
length preserving rules only in the direction indicated above.
Moreover, if we cannot apply length reducing rules to $ g = g_0\cdots g_n$
then we have $g=1$ \IFF{} both $n=0$ and $g_0 = 1$.

\section{Stallings' pregroups and their universal groups}
\label{sec:pregroup}

We now turn to the notion of pregroups in the sense of Stallings,
\cite{Stallings71}, \cite{Stallings87}.
A \emph{pregroup} $P$ is a set $P$ with a distinguished element
$\eps$,
equipped with  a partial multiplication 
$m:D \to P$, $(a,b) \mapsto ab$, where $D \subseteq P \times
P$, and an involution (or \emph{inversion}) 
$i:P \to P$, $a \mapsto a^{-1}$, satisfying the following axioms
for all $a,b,c,d \in P$. (By ``$ab$ is defined'' we mean to say that
$(a,b)\in D$ and $m(a,b)=ab$.)
\begin{enumerate}[(P1)]
\item  
$a\eps$ and $\eps a$ are defined and
$a \eps = \eps a = a;$
\item  
$a^{-1} a$ and $a a^{-1}$ are defined and $ \ \
a^{-1} a = a a^{-1} = \eps;$
\item  
if $ a b$ is defined, then so is
$b^{-1} a^{-1},$ and
$(a  b)^{-1} = b^{-1} a^{-1};$
\item  \label{it:P4}
if $a  b$ and $b
c$ are defined, then $(a b) c$ is defined if and only if $a (b
c)$ is defined, in which case
$$ (a b) c =  a (b c);$$
\item   \label{it:P5} 
if $a  b, b  c,$ and $c  d$
are all defined then either $a  b  c$ or $b
c  d$ is defined.
\end{enumerate}
It is shown in \cite{Hoare88}
that (P3) follows from (P1), (P2), and (P4), hence can be omitted.

The \emph{universal group}  $U(P)$  of the pregroup $P$ can be defined as the quotient
monoid
\begin{align*}
U(P) = 
\Gamma^*/\set{ab=c}{m(a,b)=c},
\end{align*}
where
$\Gamma = P \setminus \smallset{\eps}$ and $\eps \in P$ is identified
again with the empty word $1 \in \Gamma^*$.
The elements of $U(P)$ may
therefore represented by
finite sequences $(a_1, \ldots, a_n)$ of elements
from $\Gamma$ such that $a_ia_{i+1}$ is not defined in $P$ for $1\leq
i<n$: such sequences are called {\em $P$-reduced} sequences or {\em reduced} sequences.
Since every element in $U(P)$ has an inverse, it is clear that
$U(P)$ forms a group.

If ${\Sigma}$ is any set, then
the disjoint union $P= \smallset{\eps} \cup \Sigma \cup \ov{\Sigma}$
where $\ov{\Sigma}$ is a copy of ${\Sigma}$
yields a pregroup with involution given by $\ov\eps=\eps$, $\ov{\ov{a}}=a$,
for all $a\in \Sigma$,  such that $p \ov{p} = \eps$,
for all $p\in P$.
In this case the universal
group $U(P)$ is nothing but the free group $F(\Sigma)$.

The universal property of $U(P)$ holds trivially, namely
the canonical morphism of pregroups $P \to U(P)$ defines the
left-adjoint functor to the {forgetful}
functor from groups to pregroups.

Stallings \cite{Stallings71} showed that composition of the
inclusion map $P \rightarrow P^\ast$
 with the standard quotient map $P^\ast \rightarrow U(P)$ is injective,
where $P^\ast$ is the free monoid on $P$. The first step of his
proof establishes {\em reduced  forms} of  elements of $U(P)$, up to an equivalence
relation $\approx$ which, for completeness,  we describe here.
Define first a binary relation  $\sim$ on the set of finite sequences of elements
of $P$ by
$$(a_1,\ldots, a_i,a_{i+1},\ldots, a_n)\sim
(a_1,\ldots ,a_ic,c^{-1}a_{i+1},\ldots, a_n),$$
provided  $(a_i,c),(c^{-1},a_{i+1}) \in D.$
Then Stallings' equivalence relation $\approx$ is the transitive closure of $\sim$.

Guiding examples are again amalgamated products and HNN-extensions.
\begin{example}\label{ex:preamal}
As in Section \ref{sec:amal-prod-rewrite}, let
$A$ and $B$ be groups intersecting in a common subgroup $H$.
Consider the subset $P=A\cup B \subseteq G = A*_HB$. Define a partial
multiplication $p \cdot q$ in the obvious way;
that is $p \cdot q$ is defined if and only if either $p,q$ are both in $A$ or
$p,q$ both in $B$. Then $P$ is a pregroup  where
$D=A\times A\cup B\times B$. We
obtain the following geodesically perfect rewriting system
(where the length is computed w.r.t. $P$, thus elements of
$P$ are viewed as letters).

\[
\begin{array}{rcll}
	1  & \longrightarrow  & \eps& \;\\
p\cdot q & \longrightarrow  & r & \mbox{if } \; (p,q)\in D, pq = r \in G\\
a\cdot b & \longleftrightarrow
&ah\cdot h^{-1}b & \mbox{if } a\in A\setminus H, \; b\in B \setminus H.
\end{array}
\]
\end{example}

\begin{example}\label{ex:hnnpre}
Let $H$ be the HNN-extension HNN$(G;A,B,\Phi)$ as defined in Section~\ref{sec:hnn} and, as before, let $X$ and $Y$ be transversals for
$A$ and $B$ in $G$ with $X\cap Y = \{1\}$. Consider the subset \[
P=G\cup GtY\cup Gt^{-1}X \subset H.
\]

We define a partial multiplication by the obvious rules
(left to the reader) according to the following table.
\[
\begin{array}{rcll}
G\times G & \ras{}{} G &\\
G\times GtY & \ras{}{} GtY & \\
G\times Gt^{-1}X & \ras{}{} Gt^{-1}X & \\
GtY\times G  &\ras{}{} GtY \\
Gt^{-1}X\times G  &\ras{}{} Gt^{-1}X \\
 Gt^{-1}X \times GtY & \ras{}{} G & \mbox{if the inner
 part $XG$ is in $A$} \\
 GtY \times Gt^{-1}X & \ras{}{} G & \mbox{if the inner
 part $YG$ is in $B$}.
\end{array}
\]

This defines a pregroup $P$ for $H$, where
\[
D=G\times G
\cup G\times GtY
\cup G\times Gt^{-1}X
\cup GtY\times G
\cup Gt^{-1}X\times G
\cup S,
\]
where $S$ is the subset of $Gt^{-1}X \times GtY\cup  GtY \times Gt^{-1}X$
where inner parts $XG$ or $YG$ belong to $A$ or $B$, as appropriate.
The partial multiplication table can be directly read
from the convergent system we used in Section~\ref{sec:hnn}. As we shall
see below, it
defines an (infinite)
geodesically perfect rewriting system, where again we view
elements of $P$ as letters.
Note also that we could replace $X$ and $Y$ by $X=Y=G$ throughout the
definition of our pregroup $P$ in which the multiplication table  could
be slightly more simply described, but would be unnecessarily large.

In
\cite{Stallings71} an alternative pregroup for $H$ is defined with underlying
set consisting of equivalence classes of elements
of $G\cup t^{-1}G\cup Gt\cup t^{-1}Gt$ under the equivalence
relation generated by $t^{-1}at\sim \Phi(a)$, for $a\in A$.
However we feel that the resulting rewriting rules are obscured by the
equivalence relation on the underlying set.
\end{example}

The following is the principal result on the universal groups of pregroups.

\begin{theorem} [Stallings \cite{Stallings71}]\label{thm:sup}
Let $P$ be a pregroup. Then:
\begin{enumerate}[1)]
\item Every element of $U(P)$ can be  represented by
a $P$-reduced sequence;
\item  any two $P$-reduced sequences representing
the same element are $\approx$ equivalent, in particular they have the same length;
\item $P$ embeds into  $U(P)$.
\end{enumerate}
 \end{theorem}

\subsection{Rewriting systems for universal groups}
\label{sec:rewrite-U(P)}
The result of \cite{Stallings71} cited above may be regarded as
showing that composition of the
inclusion map $P \rightarrow P^\ast$
 with the standard quotient map $P^\ast \rightarrow U(P)$ is injective,
where $P^\ast$ is the free monoid on $P$.
We show here how to achieve this with
the help of a geodesically perfect Thue system.
Since this approach may be new we work out the details.

It is convenient to work over $P^*$ and view each
element of $P$  as a letter.
We have to distinguish whether a product is taken in  the free monoid $P^*$
or in $P$, and we introduce the following convention. Whenever we write
$[ab]$ we mean that $(a,b) \in D\subseteq P\times P$ with
$m(a,b) = [ab] \in P$: that is the
product $ab$ is defined in $P$ and yields a letter.

The system $S =S(P) \subseteq P^* \times P^*$
is now defined by the following rules.
$$
\begin{array}{rcll}
\eps & \longrightarrow  & 1 &\mbox{(= the empty word)}\\
ab & \longrightarrow  &[ab] & \mbox{if } \; (a,b) \in D  \\
ab & \longleftrightarrow  &[ac][c^{-1}b]& \mbox{if }
  \;(a,c), \, (c^{-1},b) \in D
\end{array}
$$

\begin{theorem}
\label{th:preperfect-pregroups}
Let $P$ be a pregroup. Then the following hold.
  \begin{enumerate}[1)]
 \item $P^*/S(P) \simeq U(P)$.
  \item $S$ is a geodesically perfect Thue system.
\end{enumerate}
\end{theorem}

\begin{proof} Obviously, $P^*/S$ defines  $U(P)$ which proves 1).
To prove 2) we show first that the system $S$ is strongly confluent. For this we have to consider two rules
such that the left-hand sides overlap.
Strong confluence involving only symmetric rules is trivial. Thus, we
may assume that one rule is length-reducing.
If one of the rules is $\eps \longrightarrow   1$, then (by symmetry) the
other rule is either  $\eps b \longrightarrow   b$ or
$\eps b  \longrightarrow  c[c^{-1}b]$. Since
$(c^{-1},b)\in D$ implies $(c,c^{-1}b)\in D$ and $[c(c^{-1}b)]=b$ \cite{Stallings71},
both situations lead to $b$ in at
most one step. The next situation is:
\begin{align*}
\OUT{[ab]}S{ab}{[ac]  [c^{-1}b]}
\end{align*}
Since $(a,b)$ and $(c^{-1},b)$ both belong to $D$ we
have $(a,c(c^{-1}b))\in D$, as above, and  (P\ref{it:P4})
implies that $(ac,c^{-1}b)\in D$, so
we can apply the rule
$[ac] [c^{-1}b]\longrightarrow [ab]$.
Finally, we have  to consider:
\begin{align*}
\OUT{yd}S{abd}{az}
\end{align*}
with $a,b,d \in P$ and $y, z \in P^*$.
We
may assume that one rule is length-reducing of type $ab
\longrightarrow  y=[ab]$. The other rule is 
either of type $bd  \longrightarrow  [bd]$ or
of type  $bd  \longleftrightarrow [bc][c^{-1}d]$.
Assume first that  $(b,d)\in D$, then
in both case we can use:
\begin{align*}
\IN{[ab]d}S{[abb^{-1}][bd]\,= a[bd] \, }{a[bc][c^{-1}d]}.
\end{align*}
The remaining case is that the  $(b,d)\notin D$
and  the situation is:
\begin{align*}
\OUT{[ab]d}S{abd}{a[bc][c^{-1}d]}.
\end{align*}
Since $(a,b)$, $(b,c)$ and  $(c,c^{-1}d)$ are in $D$,
(P\ref{it:P5}) implies that either
$abc$ or $bcc^{-1}d = bd$ is defined in $P$. But $bd$ is not defined,
therefore $abc$ is defined. We obtain:
\begin{align*}
\IN{[ab]d}S{[abc][c^{-1}d]}{a[bc][c^{-1}d]}.
\end{align*}

Now we show that $S$ is geodesic, from which it follows that it is geodesically perfect.  Start  with a sequence $w\in P^*$ and  apply only  length-reducing
rules until this in no longer possible. Clearly, the resulting sequence is $P$-reduced:
$w \RAS *S a_1 \cdots a_n \in \Gamma^*$ such that $a_ia_{i+1}$ is not defined in $P$ for $1\leq
i<n$. Possibly, one can still apply the symmetric rules, but we claim that any application of the symmetric rules gives again a $P$-reduced system.  Indeed, assume $ u \in \Gamma^*$ is $P$-reduced, but it is not $P$-reduced after one application of a length-preserving rule from $S(P)$.  Then there are four
consecutive elements $abde$ in $u$ and an element $c \in P$ such that
neither $ab$ nor $bd$ nor $de$ is defined, but  $bc$, $c^{-1} d $ are
defined and either   $a(bc)$ or  $(bc)(c^{-1}d)$ or
$(c^{-1}d)e$ is defined.  Assume the product $a(bc)$
is defined. Then the sequence $a, bc, c^{-1}, d$ satisfies the premise of the axiom (P5), so
either $a(bc) c^{-1}= ab $ or $(bc) c^{-1}d= bd$ must be defined,
contradicting the assumption that $u$ is $P$-reduced. Similarly,
$(c^{-1}d)e$ cannot be defined.
Suppose now that  $(bc)(c^{-1}d)$ is defined. Then the sequence $b,c,c^{-1}d$ satisfies the premise of (P4), since $(bc)$ and $c(c^{-1}d)$ are defined. Since $(bc)(c^{-1}d)$ is defined (P4) implies that  $b(c^{-1}(cd)) = b(1d) = bd$ is defined, in  contradiction with $P$-reducibility of $u$.
\end{proof}

\begin{remark}\label{stallingsoderwas}
Stallings' normal form theorem~\ref{thm:sup} 
is now a
consequence
of Theorem~\ref{th:preperfect-pregroups} because elements {}from $P$ are irreducible, and the rewriting system is geodesically perfect. Thus,
 $P$-reduced sequences that define the same elements in $U(P)$ are $\approx$ equivalent.
 \end{remark}

\begin{remark}\label{rem:strongvsnotstrong}
As above let $\Gamma = P \setminus \smallset{\eps}$.
Since $S(P) \subseteq P^* \times P^*$ is strongly confluent and geodesic, we obtain a \gp{}
presentation of the universal group $U(P)$. In some sense it is however nicer to
have such a presentation over
$\Gamma $. So, let us put
$S'(P) \subseteq \Gamma^* \times \Gamma^*$
defined by the following rules:
$$
\begin{array}{rcll}
aa^{-1} & \longrightarrow &1 & \mbox{if } \; a\in \Gamma\\
ab & \longrightarrow  &c & \mbox{if } \; (a,b)\in D,\; a\neq b^{-1},\; [ab]=c\\
ab & \longleftrightarrow  &[ac][c^{-1}b]& \mbox{if }
  \;(a,c), \, (c^{-1},b)\in D
\end{array}
$$
The difference is that a rule $aa^{-1} \longrightarrow \eps \in S$ ($\eps\in P$ is a letter)
is replaced by $aa^{-1} \longrightarrow 1 \in S'(P)$. This rule of $S'(P)$  needs
two steps of  $S(P)$, but in $S(P)$ we win strong confluence, whereas $S'(P)$ is not strongly confluent. However confluence of $S(P)$ transfers to $S'(P)$. Hence, both systems
$S(P)$ and $S'(P)$ are geodesically perfect.
\end{remark}

Using the \gp{} system $S(P)$ for $U(P)$ where we $P$ is finite we see that
the result of Rimlinger \cite{Rimlinger87b} leads to the following statement which
is slightly stronger than the result of \cite{GilHHR07}.
\begin{corollary}\label{virtfreegp}
Let $G$ be a finitely generated group. The following conditions are
equivalent.
\begin{enumerate}[1.)]
\item $G$ is virtually free.
\item $G$ can be presented by some finite \gp{} system.
\item $G$ can be presented by some finite geodesic system.
\end{enumerate}
\end{corollary}

\begin{proof}
By Rimlinger \cite{Rimlinger87b}, a finitely generated virtually free group
is the universal group $U(P)$ of some finite pregroup $P$.
By Theorem~\ref{th:preperfect-pregroups} it has a presentation by
the (finite) \gp{} system $S(P)$.
In our setting every \gp{} system is geodesic,
so we get the implication {}from 2.) to 3.) for free.
In order to pass {}from 3.) to 1.) one has to show
that the set of words which are equivalent to $1 \in G$
forms a context-free language. 
This is can be  demonstrated using an argument from \cite{die87stacs},
which has
also been used in \cite{GilHHR07}.
Consider a word $w$ and write it as as $w=uv$
such that $u$ is geodesic. The prefix  $u$ is kept on a push down stack.
Suppose that $v = a v'$, for some letter $a$.
Push $a$ onto the top of the stack: so the stack
becomes $ua$. There is no reason to suppose that $ua$ is geodesic and
we perform length reducing reduction steps on it to produce an equivalent
geodesic word $\widehat u$. Suppose this requires $k$ steps:
$$ua \RAS{k}{S_R} \widehat u$$
Let us show that we can bound $k$ by some
constant depending on only on $S$. Indeed for all letters $a$ we may fix a
word $w_a$ such that $aw_a \RAS*{S_R} 1$. But this means
$$\widehat u w_a \RAS{*}{S_R} \widetilde u,$$
where $\widetilde u$ is geodesic and $\widetilde u$ represents the same
group element as $u$ did. But $u$ was geodesic, too. Hence $\abs{u}
= \abs{\widetilde u}$. Therefore $|\widehat u|\ge |u| - |w_a|$ and this
tells us $k \leq |w_a|$. Since $k$ is bounded by some constant we see
that the whole reduction process involves a bounded suffix
of the word $ua$, only. This means we can factorise
$ua = pq$ and $\widehat u = pr$, where the length of $q$ is bounded
by some constant depending on $S$ only. Moreover,
$q \RAS{k}{S_R} r$. Since the length of $q$ is bounded this reduction
can be performed using the finite control of the
pushdown automaton. The automaton stops once  the input has been read and
then the stack gives us a geodesic corresponding to the input word $w$.
In particular, the set of words which represent $1$ 
 in the group is context-free.
 Thus, the group presented
 is context-free; and using a result of  Muller and Schupp
 \cite{ms83}
 we see that $G$ is virtually free.
\end{proof}

\subsection{Characterisation of pregroups in terms of geodesic systems}

In this section we consider Thue systems $S \subseteq \Gamma^*\times \Gamma^*$, corresponding  to group presentations, i.e., $\Gamma = X \cup X^{-1}$ and $S$ contains all the  rules $xx^{-1} \to 1, x^{-1}x \to 1, x \in X$.
We shall refer to these as {\em group rewriting systems}. 
We say that a rewriting system  $S \subseteq \Gamma^*\times \Gamma^*$  is {\em triangular} if each rule $\ell \to r \in S$ satisfies the "triangular" condition: $|\ell| = 2$, $|r| \leq 1$, so every rule in $S$ is of the form $ab \to c$ where $a, b \in \Gamma$ and $c \in \Gamma \cup \{1\}$.  Observe that a triangular system is length-reducing.

We also say that $S$ is {\em almost triangular}  if $S  = S^\prime \cup  S^\circ $, where $S^\prime$ is triangular and all rules in $S^\circ$ are trivial, i. e., of the form $a \to 1$, for some $a \in \Gamma$.
Non-trivial  examples of triangular systems come from triangulated presentations of groups. Namely, if $\langle X \mid R\rangle $ is a presentation of a group then one can {\em triangulate} this presentation by adding new generators and replacing old relations by finitely many triangular ones. 

Another type of example arises from pregroups. Let $P$ be a pregroup. In Section \ref{sec:rewrite-U(P)} we defined two rewriting systems $S(P)$ and $S^\prime(P)$ associated with $P$ that define the universal group $U(P)$. 
Notice that the length-reducing part $S^\prime(P)_R$ of $S^\prime(P)$ is triangular (here $\Gamma = P\setminus \{\eps\}$):
$$S^\prime(P)_R = \{aa^{-1}  \to 1 , ab  \to c \mid  a, b, c\in \Gamma, (a,b)\in D, [ab]=c, a\neq b^{-1}\},
$$
meanwhile, the length reducing part $S(P)_R$ of $S(P)$ is almost triangular, since it contains the trivial rule $\eps \to 1$.

Theorem \ref{th:preperfect-pregroups} implies the following result.

\begin{corollary} \label{co:triangular-pregroup}
Let $P$ be a pregroup. Then $S^\prime(P)_R$ is a triangular geodesic system,
$S(P)$ is an almost triangular geodesic system  and $U(P) = \Gamma^\ast/S^\prime(P)_R =P^\ast/S(P)_R$.
\end{corollary}
\begin{proof}
 It suffices to observe, that $S^\prime(P)_R$ and $S^\prime(P)$ define the same equivalence relation on $\Gamma^\ast$. Indeed, every rule of the type $ab \to [ac][c^{-1}b]$, where $[ac]$ and $[c^{-1}b]$ are defined, can be realized 
as a following rewriting sequence in $S^\prime(P)_R$:
 $$ab \gets  acc^{-1}b \to [ac]c^{-1}b \to [ac] [c^{-1}b],$$
 which shows that $S^\prime(P)_R$ and $S^\prime(P)$ are equivalent. The rest 
follows from Theorem \ref{th:preperfect-pregroups} and Remark 
\ref{rem:strongvsnotstrong}.
\end{proof}

To  prove the converse of this  corollary  we need some notation. Let  $S \subseteq \Gamma^\ast \times \Gamma^\ast$ be a triangular group rewriting system, where $\Gamma  = X\cup X^{-1}$. 
The congruence $\DAS{*}{S}$ on $\Gamma^\ast$ induces an equivalence
relation on the subset $\Gamma \cup \{1\}$, which we denote by $\approx$.
Define $P_S$ to be the quotient $(\Gamma\cup \{1\})/\approx$ and write
 $[z]$ for the equivalence class 
of the element $z \in \Gamma \cup \{1\}$ and in addition $\eps$ for
 the equivalence class of $1$.  
Define an involution $p \to p^{-1}$ on $P_S$ by setting $[x]^{-1} =[x^{-1}]$ 
and
$[x^{-1}]^{-1}=[x]$, for $x \in X$, and setting $\eps^{-1} = \eps$. 
(Note that, 
since $S$
is a group rewriting system, $x\approx 1$ if and only if $x^{-1}\approx 1$,
so this involution is well defined.)  
Now we define a ``partial multiplication'' on $P_S$ as follows. 
\begin{itemize}
\item
For $p, q  \in P_S \setminus \{\eps\}$   the product $pq$ is 
defined and equal to $s$ if there exist $x,y \in \Gamma$ such that 
$p=[x]$, $q=[y]$ and 
there is a rule $xy \to z \in S$,  with 
$z\in \Gamma \cup \{1\}$ 
and 
$s = [z]$. 
\item
For all $p \in P_S$ we put $p\eps =  \eps p = p$ and 
$pp^{-1} = p^{-1}p = \eps$.
\end{itemize} 
It is not hard to see that the partial multiplication on $P_S$ is 
well-defined. 

\begin{lemma}
\label{le:tringular-pregroup}
Let $S$ be a geodesic  triangular group rewriting system. Then:
\begin{itemize}
\item [1)] $P_S$ is a pregroup.
\item [2)] $U(P_S)$ is isomorphic to the group $\Gamma^\ast/S$.
\end{itemize}
\end{lemma}

\begin{proof}
Clearly, the axioms P1) and P2) hold in $P_S$ by construction. It suffices to show that P4) and P5) hold  in $P_S$, in which case  P3) follows.

Checking  P4). If any one of $p,q,r=\eps$ then P4) holds trivially,
 so we may assume that $p,q,r\in P_S\setminus\{\eps\}$.   
Suppose then $p = [a], q = [b], r = [c] \in P_S$ and the  products 
$pq$ and $qr$ are defined, i.e., $S$ contains rules 
$ab \to x$ and $bc \to y$ for some $x,y \in \Gamma \cup \{1\}$.  Suppose also that $(pq)r$ is defined in $P_S$, so either $[x] = [z]$ and $zc \to t \in S$ for some $z, t\in \Gamma \cup \{1\}$, or $[x]=\eps$, in which case let us define $t=c$. 
This means that $abc \DAS{*}{S} t$, for some $t\in \Gamma\cup \{1\}$
 and also $abc \DAS{*}{S} yc$. 
As $S$ is geodesic either 
$S$ contains a rule $yc\to u$, for
some $u\in \Gamma \cup \{1\}$,
 or 
$y=1$, in which case let us define $u=c$. 
Then $(pq)r=[t]=[u]=p(qr)$ in $P_S$.  It follows, by symmetry, that P4) holds.

P5). Again we may assume we have $p,q,r,s\in P_S\setminus \{\eps\}$ such that
 $p =[a], q= [b], r = [c], s = [d]$ and the products $pq,qr,rs$ are defined; so there are rules $ab \to x, bc \to y, cd \to z \in S$. We need to show that either $pqr$ or $qrs$ is defined. Assume $pqr$ is not defined. This means in
particular that $y\neq 1$ and that $S$ contains no rule with left 
hand side $ay$.  

We  may rewrite $abcd$ in two different ways: $abcd \to xcd \to xz$ and 
$abcd \to ayd$. As $S$ is geodesic either $S$ must contain  a rule 
which can be applied to $ayd$ or one of $a, y, d$ must be $1$. Given
our assumptions this means that $S$ contains a rule with left
hand side $yd$. Thus we have $(qr)s$ defined, so
P5) holds.   
This proves the first statement.

The second statement follows from Theorem  \ref{th:preperfect-pregroups},
Remark \ref{rem:strongvsnotstrong} and Corollary \ref{co:triangular-pregroup}. Indeed, it suffices to note that, by construction, the system $S$ is the length reducing part of the system $S^\prime(P_S)$ associated  with the pregroup $P_S$. 
\end{proof}

Combining Corollary \ref{co:triangular-pregroup} and Lemma \ref{le:tringular-pregroup} one gets the following characterisation of pregroups and their universal groups in terms of triangular geodesic systems.

\begin{theorem}
Let $P$ be a pregroup. Then the reduced part of the rewriting system $S^\prime(P)$ is a geodesic triangular group system which 
defines the universal group $U(P)$. Conversely, if $S$ is a triangular 
geodesic group system then $P_S$ is a pregroup, whose universal group is 
that defined by $S$.

\end{theorem}

This result gives a method of constructing a potentially useful pregroup for a group given by a presentation in generators and relators. It would be helpful
 to have a KB like procedure for finding such pregroups.

\begin{problem}
Design an (KB-like) algorithm that for a given finite triangular rewriting system finds an equivalent triangular geodesic system.
\end{problem}

\bibliographystyle{abbrv}
\newcommand{\Ju}{Ju}\newcommand{\Ph}{Ph}\newcommand{\Th}{Th}\newcommand{\Yu}{Y%
u}

\end{document}